\documentclass[12pt]{amsart}
\usepackage{hyperref}
\usepackage{amsfonts}
\usepackage{amsmath}
\usepackage{amssymb}
\usepackage{amscd}
\usepackage{graphics}
\usepackage{latexsym}
\usepackage{color}
\usepackage{epsfig}
\usepackage{tikz-cd}
\usepackage{amsopn}

\newtheorem{theorem}{Theorem}[section]
\newtheorem{proposition}[theorem]{Proposition}
\newtheorem{lemma}[theorem]{Lemma}
\newtheorem{corollary}[theorem]{Corollary}

\newtheorem{problem}[theorem]{Problem}

\newtheorem{conjecture}[theorem]{Conjecture}

\theoremstyle{remark}
\newtheorem{remark}[theorem]{Remark}

\newtheorem{example}[theorem]{Example}

\newtheorem{definition}[theorem]{Definition}

\newtheorem{caution}[theorem]{Caution}

\DeclareMathOperator{\cV}{\mathcal{V}}
\DeclareMathOperator{\cW}{\mathcal{W}}
\DeclareMathOperator{\cP}{\mathcal{P}}

\DeclareMathOperator{\cT}{\mathcal{T}}
\DeclareMathOperator{\cK}{\mathcal{K}}

\newcommand{\PP}{\mathbb{P}}
\newcommand{\F}{\mathbb{F}}
\newcommand{\CC}{\mathbb{C}}
\newcommand{\OO}{\mathcal{O}}

\newcommand{\ZZ}{\mathbb{Z}}

\newcommand{\RR}{\mathbb{R}}
\newcommand{\QQ}{\mathbb{Q}}

\DeclareMathOperator{\cD}{\mathcal{D}}
\DeclareMathOperator{\cA}{\mathcal{A}}

\DeclareMathOperator{\cF}{\mathcal{F}}
\DeclareMathOperator{\FF}{\mathbb{F}}
\DeclareMathOperator{\cL}{\mathcal{L}}

\def\v{\mathbf{v}}

\DeclareMathOperator{\ch}{ch}

\DeclareMathOperator{\Hom}{Hom}

\DeclareMathOperator{\Pic}{Pic}
\DeclareMathOperator{\Sym}{Sym}

\DeclareMathOperator{\rk}{rk}

\DeclareMathOperator{\Ext}{Ext}
\DeclareMathOperator{\ext}{ext}

\DeclareMathOperator{\Coh}{Coh}

\DeclareMathOperator{\Stab}{Stab}

\newcommand{\leqor}{\underset{{\scriptscriptstyle (}-{\scriptscriptstyle )}}{<}}

\begin{document}

\title{The Brill-Noether Theory of the moduli spaces of sheaves on surfaces}
\author[I. Coskun]{Izzet Coskun}
\address{Department of Mathematics, Stat. and CS \\University of Illinois at Chicago, Chicago, IL 60607}
\email{icoskun@uic.edu}
\author[J. Huizenga]{Jack Huizenga}
\address{Department of Mathematics, The Pennsylvania State University, University Park, PA 16802}
\email{huizenga@psu.edu}
\author[H. Nuer]{Howard Nuer}
\address{Faculty of Mathematics, Technion, Israel Institute of Technology, Haifa, Israel}
\email{hnuer@technion.ac.il}

\thanks{During the preparation of this article the I.C.  was partially supported by the NSF FRG grant DMS 1664296 and NSF grant DMS 2200684, J.H. was partially supported by NSF FRG grant DMS 1664303.}
\begin{abstract}
In this paper, we survey  recent developments in the Brill-Noether Theory of higher rank vector bundles on complex projective surfaces. We focus on weak Brill-Noether Theorems on rational and $K$-trivial surfaces and their applications.
\end{abstract}
\dedicatory{In honor of  Peter Newstead's 80th birthday, with great admiration}

\maketitle

\section{Introduction}

In this paper, we survey recent developments in the Brill-Noether Theory of higher rank vector bundles on complex projective surfaces. Throughout the paper we work over the field of complex numbers $\CC$. 

Brill-Noether Theory forms a cornerstone of classical curve theory. It describes representations of curves in $\PP^n$ and is a main tool in the study of moduli spaces of curves. The cohomological behavior of a general line bundle $\cL$ of degree $d$ on a smooth curve $C$ of genus $g$ is determined by $d$ and $g$ (see Proposition \ref{prop-genlinebdlcoh}). More generally, the cohomological behavior of a general stable vector bundle $\cV$ of rank $r$ and degree $d$ on $C$ is determined by $r, d$ and $g$ (see Proposition \ref{prop-curveWBNhighr}). Consequently, Brill-Noether Theory  focuses on the geometry of loci of vector bundles with unexpected cohomological behavior. 

In contrast, the cohomological behavior of a general stable bundle on a surface is not well-understood except for special surfaces or for special Chern characters. Crucially, the cohomology of a general stable bundle on a surface is not in general determined by the Euler characteristic and the slope.  Consequently, before one can study the behavior of special vector bundles, one has to understand the behavior of the generic stable bundle.  After some preliminaries \S \ref{sec-prelim} and reviewing the properties of the general stable bundles on a curve \S \ref{sec-curves}, in \S \ref{sec-pathologies} we will contrast the behavior of bundles on curves and surfaces and give examples of new phenomena that occur. 

For example, the stack of coherent sheaves of rank $r$ and degree $d$ on a curve is irreducible.  On the other hand, the stack of coherent sheaves with fixed invariants on a surface is almost never irreducible.  In many nice settings, focusing on stable sheaves results in an irreducible moduli space, although this is not always the case (see \S \ref{sec-pathologies}). Hence, stability plays a  more central role in the case of surfaces. On a curve of genus $g \geq 2$, there exist stable bundles of every rank  $r>0$ and every degree. In contrast, there are restrictions on the possible Chern characters of stable sheaves on surfaces such as the Bogomolov inequality. Furthermore, we do not have a classification of the Chern characters of stable bundles on surfaces except in special cases. Finally, the cohomology of a general line bundle with fixed invariants on a surface may not be concentrated in a single degree. We often do not know the cohomology of certain line bundles even on relatively simple surfaces such as the blowup of $\PP^2$ at 10 or more general points. This presents many challenges for studying the cohomology of  higher rank bundles.

In \S \ref{sec-WBN}, we will survey recent developments in computing the cohomology of the general stable sheaf on a surface.  We will concentrate on two techniques. First, we will describe results on rational surfaces using the stack of prioritary sheaves. Second, we will describe results on $K$-trivial surfaces using Bridgeland stability conditions. We will survey cases when weak Brill-Noether Theorems hold, i.e., when the general stable sheaf has at most one nonzero cohomology group.

 In \S \ref{sec-applications}, we will describe applications of weak Brill-Noether Theorems. We will concentrate on the following applications.
 \begin{enumerate}
 \item The classification of Chern characters of stable bundles. After recalling the Dr\'ezet-Le Potier classification of stable bundles for $\PP^2$, we will describe recent progress on Hirzebruch surfaces.
 \item The tensor product problem and the birational geometry of moduli spaces of sheaves. We will discuss the problem of computing the cohomology of the tensor product of two general stable sheaves and touch on applications to the construction of Brill-Noether divisors and the birational geometry of the moduli spaces.
 \item The classification of Ulrich bundles. We will discuss the relation between weak Brill-Noether Theorems and Ulrich bundles.
 \end{enumerate}

Finally, in \S \ref{sec-jumping}, we will survey some recent work on the cohomology jumping loci on moduli spaces of sheaves on $\PP^2$ following \cite{GouldLiu}. Very little is known (or even conjectured) about Brill-Noether loci on moduli spaces of sheaves on surfaces.  In general,   the Brill-Noether loci are reducible with components of different dimensions even in rank 1. When $c_1({\bf v})$ is a minimal positive class, however, then the Brill-Noether loci behave better as we will demonstrate.

\subsection*{Acknowledgments} We would like to thank Arend Bayer, Aaron Bertram, Lawrence Ein, Ben Gould, Joe Harris, John Kopper, Daniel Levine, Yeqin Liu,  Emanuele Macr\`i, Sayanta Mandal, Geoffrey Smith and Matthew Woolf for many stimulating discussions on the geometry of moduli spaces of sheaves.

\section{Preliminaries}\label{sec-prelim}
In this section, we collect basic facts and definitions concerning stable bundles. We refer the reader to \cite{CoskunHuizengaGokova, HuybrechtsLehn, LePotier} for more details. 

\subsection{Stability} We begin by recalling the definition of stability.
\begin{definition}
Let $X$ be a projective variety of dimension $n$ equipped with an ample  divisor $H$. The {\em $H$-slope} of a torsion free coherent sheaf $\mathcal{F}$ is defined by $$\mu_H(\mathcal{F}) = \frac{c_1 (\mathcal{F}) \cdot H^{n-1}}{ \rk(\mathcal{F})H^{n}}.$$ The sheaf $\mathcal{F}$ is called {\em $\mu_H$-(semi)stable} if for every proper subsheaf $\mathcal{E} \subset \mathcal{F}$ we have 
$$\mu_H(\mathcal{E}) \leqor \mu_H(\mathcal{F}).$$ When $H$ is fixed or irrelevant for the discussion, we will omit it from the notation.
\end{definition}

For curves, the ample $H$ does not play a role. Taking $H$ to have degree $1$, the slope becomes the degree over the rank of $\mathcal{F}$.  For higher dimensional varieties, stability depends on $H$.  Every torsion free sheaf $\mathcal{F}$ has a unique {\em Harder-Narasimhan filtration}
$$ \mathcal{F}_1 \subset \cdots \subset \mathcal{F}_n = \mathcal{F}$$ such that $\mathcal{E}_i = \mathcal{F}_i/\mathcal{F}_{i-1}$ are $\mu_H$-semistable with $$\mu_H(\mathcal{E}_1) > \mu_H(\mathcal{E}_2) > \cdots > \mu_H(\mathcal{E}_n).$$ Furthermore, $\mu_H$-semistable sheaves have a (not necessarily unique) Jordan-H\"{o}lder filtration  into $\mu_H$-stable sheaves. The reflexive hull of the associated graded sheaf is unique up to isomorphism \cite[Corollary 1.6.10]{HuybrechtsLehn}. Two semistable sheaves are called {\em $S$-equivalent} if their associated graded sheaves are isomorphic. Consequently, stable vector bundles are the building blocks of all vector bundles.

If $\cV$ and $\cW$ are two $\mu_H$-semistable sheaves and $\phi: \cV \to \cW$ is a nonzero homomorphism, then $\mu_H(\cW) \geq \mu_H(\cV)$. Furthermore, if both are stable and equality holds, then $\cV \cong \cW$. 

When $X$ is a curve,  Mumford constructed  an irreducible, projective coarse moduli space $M_{X}(r,d)$ parameterizing $S$-equivalence classes of semistable sheaves of rank $r$ and degree $d$. For higher dimensional varieties, a slightly different notion of stability is necessary. 
\begin{definition}
Let $\mathcal{F}$ be a pure dimensional coherent sheaf of dimension $d$. Then the Hilbert polynomial $P_{\mathcal{F}}(m)$ of $\mathcal{F}$ has the form 
$$P_{\mathcal{F}}(m) = \chi(\mathcal{F}(mH))= a_d \frac{m^d}{d!} + \mbox{l.o.t}.$$ The {\em reduced Hilbert polynomial} is defined by 
$p_{\mathcal{F}}(m)= \frac{P_{\mathcal{F}}(m)}{a_d}.$ The sheaf $\mathcal{F}$ is {\em $H$-Gieseker (semi)stable} if for every proper subsheaf $\mathcal{E} \subset \mathcal{F}$, we have $p_{\mathcal{E}}(m) \leqor p_{\mathcal{F}}(m)$ for $m \gg 0$.
\end{definition}
For torsion free sheaves, $a_{\dim(X)} = \rk(\mathcal{F})$. Consequently, by the Riemann-Roch Theorem, we conclude that $\mu_H$-stability implies $H$-Gieseker stability and $H$-Gieseker semistability implies $\mu_H$-semistability. When $r H^n$  and  $c_1(\mathcal{F}) \cdot H^{n-1}$ are relatively prime, then  these four notions coincide. In general, all the  implications are strict. Every torsion free sheaf has a Harder-Narasimhan filtration with respect to $H$-Gieseker semistability, and every $H$-Gieseker semistable sheaf has a Jordan-H\"{o}lder filtration with $H$-Gieseker stable quotients. Two $H$-Gieseker semistable sheaves are called {\em $S$-equivalent} if their associated graded objects are isomorphic.

Given a fixed Chern character ${\bf v}$, Gieseker \cite{Gieseker} and Maruyama \cite{Maruyama} constructed a projective coarse moduli space $M_{X, H}({\bf v})$ parameterizing $S$-equivalence classes of $H$-Gieseker semistable sheaves on $X$. Compared to the case of curves, the geometry of these moduli spaces are less well understood. The purpose of these notes is to explain aspects of the geometry of $M_{X, H}({\bf v})$ when $X$ is a surface and compare and contrast it with the case of curves.

\begin{proposition}
Let $(X, H)$ be a polarized surface. If $\cV$ and $\cW$ are $\mu_H$-semistable bundles such that $\mu_H(\cW) > \mu_H(\cV) + H \cdot K_X$, then $\Ext^2(\cV, \cW)=0$. In particular, if $H \cdot K_X <0$, then $\Ext^2(\cV,\cV)=0$. 
\end{proposition}

\begin{proof}
Let $\cV$ and $\cW$ be  $\mu_H$-semistable vector bundles. Then by Serre duality $\ext^2(\cV, \cW) = \hom (\cW, \cV(K_X))$. This is zero if $\mu_H(\cW) > \mu_H (\cV) + H \cdot K_X$ by stability.
\end{proof}

 \subsection{The Riemann-Roch Theorem} Let $X$ be a smooth, projective curve of genus $g$ and let $\cV$ be a vector bundle of rank $r$ and degree $d$. Then the Riemann-Roch Theorem computes the Euler characteristic $\chi(\cV)$
 $$\chi(\cV)= h^0(X, \cV) - h^1(X, \cV) =  d + r(1-g).$$
If $X$ is a smooth projective surface and $\mathcal{F}$ is a  sheaf  of rank $r$,   $\ch_1(\mathcal{F})=c$  and $\ch_2(\mathcal{F})=d$, then the Euler characteristic is given by 
$$\chi(\mathcal{F}) = r \chi(\OO_X) - \frac{K_X \cdot c}{2} + d,$$ where $K_X$ is the class of the canonical bundle $\omega_X$ of $X$. 

In these notes, we will be interested in $\chi(\cV \otimes \cW)$. It convenient to make a change of coordinates to express the Chern character with the logarithmic invariants. 
Given the Chern character ${\bf v}$ of a torsion free sheaf $\mathcal{F}$ on a smooth polarized surface $(X,H)$, define the total slope $\nu$ and discriminant $\Delta$ of $\mathcal{F}$ by the following formulae
$$\nu({\bf v}) = \frac{c_1({\bf v})}{\rk({\bf v})}, \quad\quad \Delta({\bf v}) = \frac{\nu({\bf v})^2}{2} - \frac{\ch_2({\bf v})}{\rk({\bf v})}.$$ Observe that these notions depend only on the Chern character of $\mathcal{F}$.  The Chern character can be easily recovered from $r, \nu$ and $\Delta$. The advantage of these invariants is that they are additive on tensor products
$$\nu(\cV \otimes \cW) = \nu(\cV) + \nu(\cW) \quad \mbox{and} \quad \Delta(\cV \otimes \cW) = \Delta(\cV) + \Delta(\cW).$$
In terms of these invariants, the Riemann-Roch formula reads
$$\chi(\cV) = r(\cV) (P(\nu(\cV)) - \Delta(\cV)),$$ where $$P(\nu)= \chi(\OO_X) + \frac{1}{2}(\nu^2 - \nu \cdot K_X)$$ is the Hilbert polynomial of $\OO_X$. 
More generally, if $\cV$ and $\cW$ are two torsion free sheaves, 
$$\chi(\cV, \cW) = \sum_{i=0}^2 (-1)^i \ext^i(\cV, \cW) = r(\cV) r(\cW)(P(\nu(\cW) - \nu(\cV))- \Delta(\cV) - \Delta(\cW)),$$ where $\ext^i(\cV,\cW) = \dim(\Ext^i(\cV,\cW))$.

\section{Brill-Noether for curves}\label{sec-curves}
In this section, a curve  will always mean a smooth, irreducible, complex projective curve. We will review  basic properties of the cohomology of stable vector bundles on curves.  In the next section, we will  contrast the cases of curves and surfaces.

\subsection{Line bundles on curves}  Let $C$ be a curve of genus $g$.  Line bundles of degree $d$ on $C$ are parameterized by $\Pic^d(C)$ which is isomorphic to a $g$-dimensional abelian variety. The Euler characteristic $\chi(\cL)=d-g+1$ determines the cohomology of the {\em general} $\cL \in \Pic^d(C)$.

\begin{proposition}[Weak Brill-Noether for line bundles]\label{prop-genlinebdlcoh}
Let $C$ be a curve of genus $g$. Let $\cL$ be a general line bundle of degree $d$ on $C$. Then 
$$h^0(C, \cL) = \max(0, d-g+1) \quad \mbox{and} \quad h^1(C, \cL) = \max(0, g-d-1).$$
\end{proposition}
\begin{proof}
If $d<0$, then $h^0(C, \cL)=0$ and $h^1(C, \cL)=g-d-1$.  If $0 \leq d \leq g-1$, let $C^{(d)}$ denote the $d$-th symmetric product of $C$. Since  $\dim(C^{(d)}) = d < \dim(\Pic^d(C))=g$, the natural map $$\phi_d: \Sym^d(C) \to \Pic^d(C)  \quad \mbox{sending} \quad  D=\sum_{i=1}^d p_i \mapsto \OO_C(D)$$ cannot be surjective. The fiber of $\phi_d$ is the linear system $|\OO_C(D)|$.  Hence, the general line bundle of degree $d$ has no global sections. 

Finally, if $d \geq g$, $h^0(C, \cL) =d-g+1 + h^1(C, \cL) \geq 1$. Letting $\omega_C$ denote the canonical bundle of $C$, by Serre duality, $h^1(C, \cL)= h^0(C, \omega_C \otimes \cL^{-1})$. By Lemma \ref{lem-genptsindependent}, general points impose independent conditions on sections of a line bundle. Hence, if $\cL= \OO_C(D)$ for a general effective divisor $D$, then $h^0(C, \omega_C(-D)) = 0$ and $h^0(C, \cL) = d-g+1$. 
\end{proof}

\begin{lemma}\label{lem-genptsindependent}
Let $X$ be an irreducible projective variety and let $\cL$ be a line bundle on $X$. If $Z$ is a general set of $m$ distinct points on $X$, then $$h^0(X, \cL\otimes I_Z) = \max(h^0(X, \cL) - m, 0).$$
\end{lemma}
\begin{proof}
If $W$ is a zero-dimensional scheme of length $m$, then $h^0(X, \cL\otimes I_W) \geq \max(h^0(X, \cL) - m, 0).$  We need to show that for a general set of points equality holds. This is trivially true if $m=0$. By induction, suppose it holds for $m \leq m_0 -1$. Since the symmetric product $X^{(m)}$ is irreducible, by the semi-continuity of cohomology, it suffices to exhibit one $Z$ of length $m_0$ for which the lemma holds. Choose $Z'$ of length $m_0 -1$ for which $$h^0(X, \cL\otimes I_{Z'}) = \max(h^0(X, \cL) - m_0+1, 0).$$ If $\cL \otimes I_{Z'}$ has no global sections, then $h^0(X, \cL\otimes I_Z)=0$ and we are done. Hence, we may assume that $\cL \otimes I_{Z'}$ has a nonzero section $s$. Take a point $p$  where $s$ does not vanish. Let $Z = Z' \cup p$. Then $h^0(X, \cL \otimes I_Z) \leq h^0(X, \cL \otimes I_{Z'}) -1$ and the proof is complete.
\end{proof}

\begin{definition}
A sheaf $\mathcal{F}$ on a projective variety $X$ satisfies {\em weak Brill-Noether} if $\mathcal{F}$ has at most one nonzero cohomology group. The sheaf $\mathcal{F}$ is {\em nonspecial} if $H^i(X, \mathcal{F})=0$ for $i >0$. Otherwise, $\mathcal{F}$ is called {\em special}.
\end{definition}

On any curve,  a general line bundle $\cL$ satisfies weak Brill-Noether and, if  $\chi(L) \geq 0$, $\cL$ is nonspecial. Consequently, classical Brill-Noether theory focuses on loci of line bundles in $\Pic^d(C)$ with unexpected cohomology. 

Recall that a  $g_d^r$ is a linear system of projective dimension $r$ and degree $d$ on the curve $C$ of genus $g$. A base-point-free $g_d^r$ corresponds to a nondegenerate morphism $C \to \PP^r$ of degree $d$. The curve $C$ admits a nondegenerate map to $\PP^r$ of degree at most $d$ if and only if $C$ has a $g_d^r$.  The Brill-Noether number is defined by  $$\rho(g,r,d) = g- (r+1)(g-d+r).$$ Let $W_d^r(C)$ be the locus of line bundles  $\cL \in \Pic^d(C)$ such that  $h^0(C, \cL) \geq r+1$. The scheme $W_d^r(C)$ has a natural determinantal structure. The celebrated Brill-Noether Theorem determines the structure of $W_d^r(C)$ on a general curve.

\begin{theorem}[Brill-Noether]\label{thm-curveBN}
Let $C$ be a curve of genus $g$ which is general in moduli.
\begin{enumerate}
\item  The curve has a $g_r^d$ if and only if $\rho(g,r,d) \geq 0$ \cite[Griffiths and Harris]{GriffithsHarris}.
\item If $\rho \geq 0$, then $W_d^r(C)$ is normal, Cohen-Macaulay of dimension $\rho$ and smooth away from $W_d^{r+1}(C)$  \cite[Gieseker]{Gieseker3}.
\item $W_d^r(C)$ is irreducible if $\rho>0$  \cite[Fulton and Lazarsfeld]{FultonLazarsfeld}.
\item When $\rho\geq0$, the universal  space  $\mathcal{W}_d^r$ of $g_d^r$'s has a unique component dominating the moduli space of curves \cite[Eisenbud and Harris]{EisenbudHarris}.
\end{enumerate}
\end{theorem}

Unlike Proposition \ref{prop-genlinebdlcoh}, the Brill-Noether Theorem requires $C$ to be general in moduli. For special curves, $W_d^r(C)$ may be reducible with components of larger than the expected dimension $\rho$.  The structure of $W_d^r(C)$ on special curves is not fully understood and is an active area of research. For example, Larson, Larson and Vogt recently described $W_d^r(C)$ for a general curve $C$ with fixed gonality \cite{LLV}.

\subsection{Higher rank vector bundles on curves} In this subsection, we will recall facts concerning stable bundles on curves and their Brill-Noether theory.

\subsubsection{Rational curves} A vector bundle $\cV$ on $\PP^1$ is isomorphic to a direct sum of line bundles $\cV \cong \bigoplus_{i=1}^r \OO_{\PP^1}(a_i) $ for a unique sequence of integers  $a_1 \leq a_2 \leq \cdots \leq a_r.$ In particular, there are no stable bundles of rank greater than one on $\PP^1$. There exists a  semistable vector bundle $\cV$ of rank $r$ and degree $d$ on $\PP^1$ if and only if $r$ divides $d$, in which case   $\cV \cong \OO_{\PP^1}(\frac{d}{r})^{\oplus r}.$ 

The cohomology of the  vector bundle $\cV = \bigoplus_{i=1}^r \OO_{\PP^1}(a_i)$  is given by 
$$h^0(\PP^1, \cV)= \sum_{a_i \geq 0} (a_i +1) \quad \mbox{and} \quad h^1(\PP^1, \cV) = \sum_{a_i \leq -2} (-1-a_i).$$

\subsubsection{Elliptic curves} Let $E$ be an elliptic curve. Vector bundles on $E$ have been classified by Atiyah \cite{Atiyah}. There exists a stable vector bundle of rank $r$ and degree $d$ on $E$ if and only if $r$ and $d$ are relatively prime. In this case, the determinant map gives an isomorphism between $M_E(r,d)$ and $E$. In particular, by taking direct sums, there exist semistable bundles of every rank $r$ and degree $d$ on $E$.  

If $\cV$ and $\cW$ are two semistable bundles on $E$ with $\mu(\cV) > \mu(\cW)$, then $$\Ext^1(\cW, \cV) \cong \Hom(\cV, \cW)^*=0$$ by semistability. By induction on the length of the Harder-Narasimhan filtration, we conclude that a vector bundle on $E$ is a direct sum of  semistable vector bundles. Hence, the cohomology of any bundle on $E$ is determined by the cohomology of the bundles appearing in its Harder-Narasimhan filtration.

\begin{proposition}\label{prop-genus1}
Let $\cV$ be a semistable bundle on an elliptic curve $E$. Then $\cV$ satisfies weak Brill-Noether if and only if  $\OO_E$ is not a Jordan-H\"{o}lder factor of $\cV$.
\end{proposition}
\begin{proof}
Let $\cV$ be a stable vector bundle on $E$. If $\mu(\cV)<0$, then $H^0(E,\cV)=0$ by stability. Similarly, if $\mu(\cV)>0$, then $H^1(E,\cV)=H^0(E, \cV^*)=0.$ If $\mu(\cV)=0$ and $\cV$ is stable, then $\cV$ must be a line bundle of degree $0$. The cohomology of $\cV$ in that case vanishes except when $\cV= \OO_E$. We conclude that on an elliptic curve all stable bundles except $\OO_E$ satisfy weak Brill-Noether. By considering the Jordan-H\"{o}lder filtration, we conclude that any semistable bundle $\cV$ on $E$ satisfies weak Brill-Noether provided that $\OO_E$ is not a Jordan-H\"{o}lder factor of $\cV$. It is easy to see that if $\OO_E$ is a Jordan-H\"{o}lder factor, then $\cV$ must have nonvanishing cohomology.
\end{proof}

\begin{example}\label{ex-elliptic}
The cohomology of a semistable bundle whose Jordan-H\"{o}lder factors contain $\OO_E$ depends on the extension class. Let $$0 \to \OO_E \to \cV \to \OO_E \to 0$$ be the unique nontrivial extension of $\OO_E$ by $\OO_E$. Then $h^0(E, \cV)= h^1(E, \cV)=1$, whereas  $h^0(E, \OO_E \oplus \OO_E) = h^1(E, \OO_E \oplus \OO_E)=2$. 
\end{example}

 \subsubsection{Curves of higher genus} For the rest of this section, let  $C$ be a curve of genus $g \geq 2$. The following theorem describes the moduli spaces of semistable sheaves on $C$.
 
\begin{theorem}\label{thm-curveirreddim}
Let $C$ be a smooth curve of genus $g\geq 2$. 
\begin{enumerate}
\item Then the stack of coherent sheaves $\mathcal{C}oh_{d,r}$ of degree $d$ and rank $r$ on $C$ is a smooth, irreducible Artin stack.
\item  The moduli space  $M_C(r,d)$ is an irreducible, projective variety of dimension $r^2(g-1) + 1$, which contains the locus of stable bundles as a dense open set.  It is smooth at points corresponding to stable bundles. In particular, if $r$ and $d$ are relatively prime, then $M_C(r,d)$ is smooth. 
\end{enumerate} 
\end{theorem}

\begin{proof}[Sketch of proof]
We refer the reader to \cite{Hoffmann} for a detailed proof of (1). Briefly, using the Grothendieck Quot scheme, one can show that $\mathcal{C}oh_{d,r}$ is algebraic. Obstructions to deformations of a coherent sheaf $\mathcal{F}$ are contained in $\Ext^2(\mathcal{F},\mathcal{F})$. Since  $\Ext^2(\mathcal{F},\mathcal{F})=0$ on $C$,  $\mathcal{C}oh_{d,r}$ is smooth of dimension $-\chi(\mathcal{F}, \mathcal{F})$. Finally, by induction on rank and degree, one can show that $\mathcal{C}oh_{d,r}$ is connected to conclude (1).  The irreducibility of $M_C(r,d)$ (when nonempty) follows from openness of semistability.

We now construct semistable vector bundles of rank $r$ and degree $d$ on curves of genus $g \geq 1$. First, assume $r$ and $d$ are coprime. 
Let $\frac{e}{s} < \frac{d}{r} < \frac{f}{t}$ be the three consecutive fractions in the Farey sequence for $r$. This implies that the middle fraction is the \emph{mediant} of its two neighbors: $$\frac{d}{r} = \frac{e+f}{s+t}.$$ Then $s$ and $t$ are strictly less than $r$ and the $e,s$ (similarly $f,t$) are  coprime. By induction on the rank, there exist stable bundles $\cV_1$ and $\cV_2$ of rank and degree $(s,e)$ and $(t,f)$, respectively. Since $$\chi(\cV_2, \cV_1) = et-sf + r(1-g) < 0,$$ there exists a nonsplit extension of the form $$0 \to \cV_1 \to \cV \to \cV_2 \to 0.$$ We claim that $\cV$ is stable. Otherwise, let $\cW$ be the maximal destabilizing subbundle. Since $\mu(\cW) > \mu(\cV)$, $r(\cW) < r(\cV)$ and $\frac{f}{t}$ is the next smallest fraction in the Farey series for $r$, we must have $\mu(\cW) \geq \mu(\cV_2)$. Hence, either $\mu(\cW) > \mu(\cV_2)$ and the natural  map from $\cW$ to $\cV_2$ is zero or $\mu(\cW)= \mu(\cV_2)$. In the first case, there is an induced map from $\cW$ to $\cV_1$, which is a contradiction since $\cW$ is stable and $\mu(\cV_1) < \mu(\cW)$.  In the second case, we must have $\cW = \cV_2$ and we get a splitting of the sequence, contrary to our assumptions.

Hence, when $r$ and $d$ are  coprime, the moduli space $M_C(r,d)$ is nonempty. Furthermore, if $\cV$ is a stable bundle, then $$\hom(\cV,\cV)=1, \quad \mbox{and} \quad \ext^1(\cV,\cV)= r^2(g-1) + 1.$$ By basic deformation theory, $T_{[\cV]} M_C(r,d)= \Ext^1(\cV,\cV)$ and $M_C(r,d)$ is smooth at $\cV$ of dimension $r^2(g-1)+1$.

If $r$ and $d$ are not  coprime, let $(r, d) = (kr', kd')$, where $r'$ and $d'$ are  coprime. There exists a stable vector bundle $\cV$ of rank and degree $(r', d')$.  Taking the direct sum of  $k$ copies of $\cV$, we obtain a semistable bundle of rank $r$ and $d$. Hence, if $g \geq 1$, $M_C(r,d)$ is nonempty. 

Now let  $g \geq 2$. By induction on the gcd $k$, assume that there are stable bundles for all $k' < k$. The dimension of the stack $\mathcal{M}_C(r,d)$ is $r^2(g-1)$. Consider the first step of a Jordan-H\"{o}lder filtration
$$0 \to \cV_1 \to \cV \to \cV_2 \to 0$$ with the rank and degree of $\cV_i$ equal to $(k_i r', k_i d')$ and $k_1 + k_2 = k$. We have  $$\ext^1(\cV_2, \cV_1) = k_1 k_2 (r')^2 (g-1) + \hom(\cV_2, \cV_1).$$
Since $\cV_1$ is stable, either $\hom(\cV_2, \cV_1)=0$ or $\cV_1$ is a factor of the  Jordan-H\"{o}lder filtration of $\cV_2$. In the first case, the dimension of such extensions is bounded above by
$$(k_1r')^2 (g-1) + k_1k_2 (r')^2 (g-1) + (k_2r')^2 (g-1) < r^2(g-1).$$ In the latter case, $\cV_1$ is determined up to finitely many choices and $\hom(\cV_2, \cV_1)\leq \frac{k_2}{k_1}.$ Hence such loci is bounded by $$\frac{k_2}{k_1} + k_1k_2 (r')^2 (g-1) +(k_2r')^2 (g-1) < r^2(g-1).$$ We conclude that there must be stable bundles of rank $r$ and degree $d$ on $C$.
\end{proof}

\begin{proposition}[Weak Brill-Noether for curves]\label{prop-curveWBNhighr}
Let $\cV$ be a general semistable vector bundle of rank $r$ and degree $d$ on a curve $C$ of genus $g$. Then $\cV$ satisfies weak Brill-Noether:
$$h^0(\cV) = \max(0, d-r(g-1)), \quad \mbox{and} \quad h^1(\cV) = \max(0, r(g-1)-d).$$
\end{proposition}

\begin{proof}
The proposition is easy for  $g=0$ and follows from Proposition \ref{prop-genus1} for $g=1$. Hence, we may assume that $g\geq 2$.  By the semicontinuity of cohomology and Theorem \ref{thm-curveirreddim}, it suffices to exhibit one coherent sheaf of rank $r$ and degree $d$ with the expected cohomology. By Riemann-Roch, the Euler characteristic is given by $\chi(\cV) = d + r(1-g)$. 
By the division algorithm, write $d=qr+s$ with $0 \leq s < r$. Consider a direct sum of line bundles
$$\cV = \bigoplus_{i=1}^{r-s} \cL_i\oplus  \bigoplus_{j=1}^s \cL_j', $$ where $\cL_i$ are general line bundles of degree $q$ and $\cL_j'$ are general line bundles of degree $q+1$. By Proposition \ref{prop-genlinebdlcoh}, $h^0(\cL_i) = \max(0, q-g+1)$ and $h^0(\cL_j') = \max(0, q+2-g)$. Since $h^0(\cV) = \sum_{i} h^0(\cL_i) + \sum_j h^0(\cL_j')$, we conclude that $h^0(\cV)=0$ if $q+2 \leq g$ or if $s=0$ and $q+1 \leq g$.  Finally, we have $h^0(\cV) = d + r(1-g)$ if $q+1 \geq g$. This proves the statement for $h^0$. The statement for $h^1$ follows by Serre duality.
\end{proof}

\begin{caution}
We warn the reader that points of $M_C(r,d)$ correspond to $S$-equivalence classes of bundles. Certain properties of bundles such as cohomology or global generation are not invariant under $S$-equivalence. The two bundles in Example \ref{ex-elliptic} are $S$-equivalent, but they have different cohomology groups. Similarly, $\OO_E \oplus \OO_E$ is globally generated, whereas the nontrivial extension is not.  
\end{caution}

Since the cohomological behavior of the general stable  bundle is as expected, Brill-Noether theory concentrates on the cohomology jumping loci.  In higher rank, much less is known about cohomology jumping loci. We refer the reader to \cite{GM} for a  survey. As in the case of line bundles, there are no interesting jumping loci as soon as the slope of the vector bundle is sufficiently large.

\begin{proposition}
Let $C$ be a curve of genus $g$. If  $\mu > 2g-2$, then every bundle $\cV$ in $M_C(r,d)$ is nonspecial.
\end{proposition}

\begin{proof}
By Serre duality, $$h^1(C, \cV) = h^0(C, \cV^*\otimes \omega_C) = \hom(\cV, \omega_C).$$ Since $\cV$ and $\omega_C$ are semistable and $$\mu(\cV) > 2g-2 = \mu(\omega_C),$$ we conclude that $\hom(\cV, \omega_C)=0$. Hence, every semistable bundle of slope greater than $2g-2$ has no higher cohomology.
\end{proof}

Consequently, the Brill-Noether problem is only interesting when $\mu(V) \leq 2g-2$. Since $h^1(C, \omega_C) =1$, the bound in the proposition is sharp.

Several other properties of bundles such as global generation and ampleness play a central role in geometry. Recall that a vector bundle $\cV$ is called \emph{ample} if $\OO_{\PP \cV}(1)$ is an ample line bundle on $\PP \cV$, or equivalently, if for every coherent sheaf $\cF$ there exists an integer $n_0$ such that $\cF\otimes \Sym^n(\cV)$ is globally generated for all $n\geq n_0$. For curves, it is easy to determine when the general stable bundle has these properties.

\begin{proposition}
Let $C$ be a  curve of genus $g \geq 2$.
\begin{enumerate}
\item \cite{Laumon, Sundaram} A general stable bundle of rank $r$ and degree $d$ on $C$ is globally generated if $d  - rg \geq 1$.
\item If $\mu > 2g-1$, then every semistable vector bundle of rank $r$ and degree $d$ on $C$ is globally generated.
\item \cite[Theorem 2.4]{Hartshorneample} If $\mu >0$, then every semistable bundle of rank $r$ and degree $d$ on $C$ is ample.
\end{enumerate}
\end{proposition}

\begin{remark}
Global generation is not in general an open property. For example, let $\cL$ be a line bundle of degree $0$ on a positive genus curve $C$. Then $\cL$ is globally generated if and only if $\cL \cong \OO_C$. Hence, the locus of globally generated degree 0 line bundles is a single point in a $g$-dimensional abelian variety. However,  the locus of  globally generated bundles  is open in the locus of nonspecial bundles. 
\end{remark}

\begin{proof}
If $d-rg \geq 1$, then by Proposition \ref{prop-curveWBNhighr},  the general stable bundle has no higher cohomology and has at least $r+1$ sections. By the irreducibility of $\Coh_{d,r}$, it suffices to exhibit one globally generated nonspecial  bundle of  rank $r$ and degree $d$.  

A general line bundle $\cL$ of degree $d \geq g+1$ has the property that $h^1(C, \cL(-p))=0$ for every $p \in C$. Indeed, consider the locus 
$$X_{d-1} := \{ D \in C^{(d-1)} | h^0(C, D) \geq d-g + 1 \} \subset C^{(d)}.$$ By Proposition \ref{prop-curveWBNhighr}, the codimension of $X_{d-1}$ is at least $1$. The fibers of the natural map $\phi_d: C^{(d-1)} \to \Pic^{d-1}(C)$ of $X_{d-1}$ have  dimension at least $d-g$.  Consequently, the image of $X_{d-1}$  has dimension at most $g-2$.  Hence, the codimension of the locus of line bundles of degree $d-1$ with unexpected cohomology is at least $2$. Therefore, for the general $\cL \in \Pic^{d}(C)$, $h^0(C, \cL(-p))= h^0(C, \cL) -1$ for every $p \in C$. By Serre duality, $h^1(C,\cL(-p))=0$ for every $p \in C$.

By the Euclidean algorithm, write $d= rm + s$. If $m \geq g+1$, then we can take $\cV$ to be an extension of the form 
$$0 \to \cL^{r-s} \to \cV \to \cK^{s} \to 0,$$ where $\cL$ and $\cK$ are general line bundles of degree $m$ and $m+1$, respectively. We then have that $\cV$ is nonspecial and $H^1(C, \cV(-p))=0$ for every $p \in C$. Taking the long exact sequence associated to $$0 \to \cV(-p) \to \cV \to \cV|_p \to 0,$$ we conclude that $H^0(C, \cV)$ surjects onto $H^0(C, \cV|_p)$ for every $p \in C$ and $\cV$ is globally generated. If $m=g$, a more subtle argument is needed. We claim that the Brill-Noether locus in $M_C(r,d)$ parameterizing special vector bundles has codimension at least $2$ if $d - r (g-1) \geq 1$ (see  \cite{Laumon, Sundaram}). Granting this claim, $\deg(\cV(-p)) - r(g-1) \geq 1$. Consequently, $\cV(-p)$ does not have $h^1$ for the general $\cV$ and every $p \in C$. Hence, $\cV$ is a globally generated,  nonspecial bundle.  To prove the claim, recall that 
the fibers of the determinant map $M_{C}(r, d) \to \Pic^d(C)$ are Fano varieties of Picard rank 1 \cite{DrezetNarasimhan}. Hence, it suffices to exhibit a complete curve in each fiber which is nonspecial. Furthermore, by taking direct sums, it suffices to do this when the rank and the degree are coprime. Now it is easy to construct such complete families by induction on the rank as  in the proof of Theorem \ref{thm-curveirreddim}.

Now assume $\cV$ is semistable and $\mu(\cV)> 2g-1$. For any point $p \in C$, consider the exact sequence
$$0 \to \cV(-p) \to \cV \to \cV_p \to 0.$$ We have that $H^1(C, \cV(-p)) \cong H^0(C, \cV^*\otimes \omega_C(p))$. Since $\mu(\cV^*\otimes \omega_C (p)) < 0$ and is semistable, it cannot have any sections. Hence the evaluation map $H^0(C, \cV) \to H^0(\cV|_p)$ is surjective for every point $p$ and $\cV$ is globally generated. 

Hartshorne \cite{Hartshorneample} proves (3) by showing that a vector bundle $\cV$ of positive degree on a curve all of whose symmetric powers $\Sym^n(\cV)$ are semistable is ample. In characteristic 0, by the Narasimhan-Seshadri Theorem, the symmetric powers of a semistable bundle are semistable. More generally, Hartshorne uses this fact to show that any vector bundle of positive degree on a smooth curve all of whose quotients have positive degree is ample.
\end{proof}

The locus of stable vector bundles that fail to be globally generated has been studied by Kopper and Mandal \cite{KopperMandal}. Their main theorem is the following.
\begin{theorem}\cite[Theorem 1.2]{KopperMandal}
 Let $N_C(r,d)$ be the locus of {\em stable} vector bundles of rank $r \geq 2$ and degree $d$ on a smooth curve of genus $g \geq 2$ that fail to be globally generated. Assume $rg+1\leq d \leq r(2g-1)-1$.
 \begin{enumerate}
 \item $N_C(r,d)$ is nonempty and it has a component of expected codimension $d-rg$ and no component of smaller codimension.
 \item If $rg+g-1 \leq d \leq r(2g-1)-1$, then $N_C(r,d)$ is irreducible of the expected codimension $d-rg$.
 \end{enumerate}
\end{theorem}

\section{Pathologies for surfaces}\label{sec-pathologies}

In this section, we will contrast the behavior of stable bundles on surfaces with that of curves discussed in the last section. We will see that many of the nice properties that are valid for curves fail for surfaces. This makes studying Brill-Noether theory on surfaces much more challenging. Here are some of the new difficulties that arise.

\begin{enumerate}
\item A general line bundle with fixed invariants on a surface can have more than one nonzero cohomology group. Hence, unlike on curves, the cohomology of a general line bundle on a surface is more delicate and can be challenging to compute.
\item Unlike in the case of curves of genus $g\geq2$, there are restrictions on the possible invariants of stable sheaves on surfaces. Determining the Chern characters of stable sheaves on a surface is a challenging problem. 
\item Unlike in the case of curves, the stack of coherent sheaves with fixed invariants on a surface is almost always reducible. Hence, stability plays a much greater role in determining the cohomology of a `general' stable sheaf on a surface. Furthermore, on a surface $X$, the moduli space $M_{X, H}({\bf v})$ itself may be reducible and even disconnected. Hence, it does not always even make sense to talk about a general stable sheaf with fixed invariants. 
\end{enumerate}

In this section, we will give some simple examples of these phenomena.

\subsection{The cohomology of line bundles on surfaces}   Proposition \ref{prop-genlinebdlcoh} showed that a general line bundle of degree $d$ on  {\em any} smooth curve has at most one nonzero cohomology group. This fails for surfaces. The general line bundle with a given class on a surface may have two or three nonzero cohomology groups. Moreover, whether weak Brill-Noether holds for a line bundle may depend on the isomorphism class of the surface.

\begin{example} We give several easy examples where weak Brill-Noether fails for line bundles on surfaces.
\begin{enumerate}
\item Let $X$ be a K3 surface. Then the Picard group of $X$ is discrete and $$h^0(X, \OO_X)=h^2(X, \OO_X)=1.$$
\item Let $X$ be a very general surface of degree $d \geq 4$ in $\PP^3$. By the Noether-Lefschetz Theorem, the Picard group of $X$ is generated by $\OO_X(1)$. Let $0 \leq k \leq d-4$. The long exact sequence associated to 
$$0\to \OO_{\PP^3}(k-d) \to \OO_{\PP^3}(k) \to \OO_X(k) \to 0$$ implies that 
$$h^0(X, \OO_X(k)) = \binom{k+3}{3}, \quad h^1(X, \OO_X(k))=0, \quad  h^2(X, \OO_X(k))= \binom{d-k-1}{3}.$$ Observe that both $h^0$ and $h^2$ are nonzero.
\item Let $Y \subset \PP^3$ be a very general surface of degree $d \geq 5$ and $\ell$ be a general line in $\PP^3$.  Let $X$ denote the blowup of $Y$ along $\ell \cap Y= \{ p_1, \dots, p_d\}$ and  let $E_i$ denote the exceptional divisor over $p_i$.  Set $E = \sum_{i=1}^d E_i$. Let $H$ be the pullback of the hyperplane class from $Y$.  Then the Picard group of $X$ is discrete.  For an integer $1 \leq k \leq d-4$, the sections of $\OO_X(kH-E)$ are given by hypersurfaces of degree $k$ that vanish along $\ell \cap Y$. By B\'ezout's Theorem, these hypersurfaces of degree $k$ must vanish along $\ell$. By Serre duality, $$h^2(X, kH-E)= h^0(X, (d-4-k)H + 2E).$$ Since $E_i \cdot E_j=-\delta_{i,j}$, $2E$ must be in the base locus of $\OO_X((d-4-k)H + 2E)$. Consequently, $h^2(X, kH-E)= h^0(X, (d-4-k)H)$.  By an easy Euler characteristic computation, we conclude that 
$$ h^0(X, kH-E) = \binom{k+3}{3} - k -1, \ h^1(X, kH - E) = d-k-1,$$  $$h^2(X, kH-E) = \binom{d-k-1}{3}.$$ Hence, all three cohomology groups are nonzero.
Observe that if $X$ were the blowup of $Y$ along $d$ general points, the cohomology would instead be
$$ h^0(X, kH-E) = \max\left(0, \binom{k+3}{3} - d\right),  \ h^1(X, kH - E) = \max\left(0, d-\binom{k+3}{3}\right),$$  $$h^2(X, kH-E) = \binom{d-k-1}{3}.$$ Hence, the generic cohomology of a line bundle depends on the isomorphism class of the surface.
\item Let $A$ be an abelian surface of Picard rank $1$. Let $X$ be the blowup of $A$ at a point with exceptional divisor  $E$. Let $H$  denote the pullback of a very ample divisor on $A$. Let $m$ be a nonnegative integer. Then any line bundle algebraically equivalent to $\OO_X(H + mE)$ is of the form $\mathcal{L} \otimes \OO_X(H + mE)$, where $\mathcal{L}$ is the pullback of a line bundle of degree $0$ from $A$.  The line bundle $\mathcal{L} \otimes \OO_X(H)$ is the pullback of an ample line bundle from $A$, which has positive Euler characteristic and no higher cohomology by the Kodaira Vanishing Theorem.  Consequently, $$h^0(X, \mathcal{L} \otimes \OO_X(H+mE)) >0.$$ However, $$\chi(\mathcal{L} \otimes \OO_X(H+mE)) = \chi(\OO_X) + \frac{1}{2} (H^2 - m^2 -m) <0$$ if $m \gg 0$. We conclude that these line bundles have nontrivial $h^0$ and $h^1$. 
\end{enumerate}
\end{example}

In general, computing  the cohomology of a line bundle on a surface can be a  hard problem. For example, already for the blow-up of $\PP^2$ in $m \geq 10$ very general points, we do not know the dimensions of the spaces of global sections of all line bundles. In this case there is a precise conjecture due Segre, Harbourne, Gimigliano and Hirschowitz. 

\begin{conjecture}[SHGH Conjecture]\label{conj-SHGH}\cite{Segre, Harbourne, Gimigliano, Hirschowitz}
Let $X$ be a very general blow up of $\PP^2$ at $m \geq 10$ points $p_1, \dots, p_m$. Let $H$ denote the pullback of the class of a line and let $E_i$ denote the exceptional divisor over $p_i$. Let $D= dH - \sum_{i=1}^m n_i E_i$ be a divisor with $d > 0$ and $n_i \geq 0$. Then $$h^0(X, \OO_X(D)) = \binom{d+2}{2} - \sum_{i=1}^m \binom{n_i+1}{2}$$ if and only if $\OO_X(D)$ does not have a multiple $(-1)$-curve in its base locus. 
\end{conjecture}
The SHGH Conjecture, if true, provides an efficient algorithm for computing the cohomology of any line bundle on $X$. Despite steady progress, the conjecture remains open in general.

Given that weak Brill-Noether does not hold in general for line bundles on surfaces, we do not expect it to hold for higher rank bundles either. Unlike the case of curves, for surfaces already computing the cohomology of a general stable sheaf in a component of the moduli space is an interesting and challenging problem.  

\begin{problem}[Weak Brill-Noether]
Given an irreducible component of $M_{X, H}({\bf v})$, compute the cohomology of the general sheaf $\cF$ in that component.
\end{problem}

\subsection{The cohomology of rank 1 sheaves} Once one knows the cohomology of line bundles on a surface $X$, understanding the general cohomology of rank 1 sheaves does not present new difficulties.  If $\cV$ is a rank 1 torsion-free coherent sheaf on a smooth projective surface, then $\cV \cong I_Z \otimes \mathcal{L}$, where $\mathcal{L}$ is a line bundle and $I_Z$ is the ideal sheaf of a zero-dimensional subscheme of $X$. The Hilbert scheme $X^{[n]}$ parameterizing length $n$ zero-dimensional subschemes of $X$ is a smooth, irreducible, projective variety of dimension $2n$.  The long exact sequence associated  to the standard exact sequence
$$0 \to I_Z \otimes \cL \to \cL \to \OO_Z \to 0$$ implies that $H^2(X, I_Z \otimes \cL ) \cong H^2(X, \cL)$.
If $Z$ is a general set of $n$ points on $X$, then the map $H^0(X, \cL) \to H^0(X, \OO_Z)$ has maximal rank. We conclude that $H^0(X, I_Z \otimes \cL) =0$ for a general set of points if and only if $n \geq H^0(X, \cL)$. Similarly, $h^1(X, I_Z \otimes \cL) = h^1(X, \cL) +n- h^0(X,\cL)$ for a general set of points $Z$. Hence, we can compute the cohomology of $I_Z \otimes \cL$ for a general set of points $Z$ purely based on the cohomology of $\cL$.

\begin{remark}
We caution the reader that rank 1 sheaves behave slightly differently from higher rank sheaves. When $n>0$, the moduli space consists entirely of non-locally-free sheaves. Consequently, one cannot directly apply Serre duality. In fact, if $\cL$ is a line bundle with $h^2(X, \cL) \not= 0$, then $\cL\otimes I_Z$ has nonvanishing $h^1$ and $h^2$ as soon as $n>h^0(X, \cL)$. In particular, if $\cL$ is an  ample line bundle on $X$ with $h^0(X, \cL) \not=0$, then $I_Z \otimes \omega_X \otimes \cL^*$ has nonvanishing $h^1$ and $h^2$  as soon as $n \geq 1$. In particular, one can only hope for weak Brill-Noether to hold  for rank 1 sheaves when $h^2(X, \cL) \not= 0$.

When the rank $r >1$, there may be components of $M_{X, H}({\bf v})$ that consist entirely of non-locally-free sheaves. However, such components do not exist for certain surfaces such as $\PP^2$ and Hirzebruch surfaces. Moreover, when $\Delta({\bf v}) \gg 0$, then the general member of the moduli space parameterizes a vector bundle. Consequently, one may apply Serre duality to compute the cohomology of the general sheaf.
\end{remark}

\subsection{The moduli space may be empty}  Theorem \ref{thm-curveirreddim} shows that on a curve of genus $g \geq 2$, the moduli space $M_C(r,d)$ is nonempty, irreducible and of the expected dimension for every $r \geq 1$ and every $d$. In contrast, there are restrictions on the possible Chern characters of stable sheaves on surfaces. For example, the Bogomolov inequality provides an important constraint.

\begin{theorem}[The Bogomolov inequality]
Let $(X, H)$ be a polarized smooth surface and let $\mathcal{F}$ a $\mu_H$-semistable sheaf. Then $\Delta(\mathcal{F}) \geq 0$.  
\end{theorem}

For abelian surfaces Bogomolov's inequality completely characterizes stable Chern characters. For some surfaces such as $\PP^2$ or K3 surfaces one can easily strengthen the Bogomolov inequality. 
\begin{example} Let  $X= \PP^2$ and let $\cV$ be a torsion-free stable sheaf. Then $$\ext^2(\cV,\cV)=\hom(\cV, \cV(-3))=0$$ by Serre duality and stability. Therefore, $$\chi(\cV,\cV) = r^2 (1- 2\Delta(\cV))= \hom(\cV,\cV)- \ext^1(\cV,\cV) \leq 1.$$ Consequently, $$\Delta (\cV) \geq \frac{1}{2} \left(1 - \frac{1}{r^2}\right).$$ We will later see that there are stronger restrictions on the Chern characters of stable sheaves on $\PP^2$.
\end{example}

 \begin{example}
 Let $X$ be a K3 surface and let $\cV$  be a torsion-free stable sheaf. Since  $\hom(\cV,\cV)= \ext^2(\cV,\cV) =1$, we have
$$\chi(\cV,\cV) = r^2 (2- 2\Delta(\cV)) \leq 2.$$
Hence, $\Delta(\cV) \geq 1 - \frac{1}{r^2}$. On K3 surfaces this inequality characterizes the Chern characters of stable sheaves.
\end{example}

We will say that a Chern character is $\mu_H$-(semi)stable  (respectively, $H$-Gieseker (semi)stable) if there exists a $\mu_H$-(semi)stable  (respectively, $H$-Gieseker (semi)stable) sheaf with that Chern character. The following is a central  open problem.

\begin{problem}\label{prob-classifystabch}
Given a polarized surface $(X,H)$ classify the Chern characters of $\mu_H$ or $H$-Gieseker semistable sheaves. Determine when there exist stable sheaves with the given Chern character.
\end{problem}

The results concerning Problem \ref{prob-classifystabch} have two flavors. First, there are important asymptotic results for all surfaces guaranteeing the existence of stable bundles with sufficiently large discriminant. 

\begin{theorem}[O'Grady]\label{thm-OGrady}
Let $(X,H)$ be a smooth polarized surface and let ${\bf v}$ be a Chern character with $r({\bf v}) >0$. If $\Delta_H({\bf v}) \gg 0$ (where the necessary inequality depends on $r$, $X$ and $H$), then the moduli space $M_H({\bf v})$ is normal, generically smooth, irreducible and nonempty of the expected dimension. Furthermore, the slope stable sheaves are dense in $M_H({\bf v})$.
\end{theorem}

Second, Problem \ref{prob-classifystabch} has been studied and solved on certain surfaces including $\PP^2$ \cite{DrezetLePotier}, K3 surfaces \cite{Yoshioka99}, Abelian surfaces \cite{Yoshioka23}, Hirzebruch surfaces \cite{CoskunHuizengaHExist}, Enriques surfaces \cite{Nuer16a, Nuer16b, NuerYoshioka}, elliptic and bielliptic surfaces. There has been some progress for certain other classes of surfaces such as quintic and sextic surfaces in $\PP^3$ at least in rank 2 \cite{MestranoSimpson}. However, the general problem remains wide open.  

Not knowing the Chern characters of stable sheaves often presents difficulties in constructions. One is often forced to assume inequalities on the discriminant in constructions. There is however a general way to construct rank 2 bundles via the Serre construction, which we now recall.

 \subsection{The Serre Construction}  Let $X$ be a smooth projective surface. Let $Z\subset X$ be a zero-dimensional local complete intersection scheme of length $n$. Then $Z$ satisfies the {\em Cayley-Bacharach property} with respect to a line bundle $\mathcal{L}$ on $X$ if for any subscheme $W \subset Z$ of length $n-1$, any section of $\mathcal{L}$ vanishing on $W$ vanishes on $Z$. The name is inspired by the classical Cayley-Bacharach Theorem which asserts that any cubic curve which contains 8 of the 9 intersection points of two cubic curves in $\PP^2$  also contains the ninth.
 
 \begin{theorem}[The Serre correspondence]\label{thm-Serre}
 Let $X$ be a smooth projective surface and $Z \subset X$ be a local complete intersection subscheme of dimension zero and length $n$. Let $\mathcal{L}$ be a line bundle on $X$. Then there exists an extension 
 $$0\to \OO_X \to \mathcal{V} \to I_Z \otimes \mathcal{L} \to 0$$ with $\cV$ locally free if and only if $Z$ satisfies the Cayley-Bacharach property for $\omega_X \otimes \mathcal{L}$.
 \end{theorem}
 
 The Serre construction allows one to  construct stable vector bundles of rank 2 on surfaces.

 \begin{corollary}
 Let $H$ be an ample divisor on $X$. Let $\mathcal{L}$ be a line bundle such that $c_1(\mathcal{L})\cdot H > 0$. Let $Z$ be a nonempty set of distinct points such that $Z$ satisfies the Cayley-Bacharach property for $\omega_X \otimes \mathcal{L}$. Assume that $H^0(X, \mathcal{N} \otimes I_Z) =0$ for any line bundle $\mathcal{N}$ with $c_1(\mathcal{N}) \cdot H \leq c_1(\mathcal{L}) \cdot H$. Then the general extension $\cV$
 $$0\to \OO_X \to \mathcal{V} \to I_Z \otimes \mathcal{L} \to 0$$ is a $\mu_H$-stable bundle.
 \end{corollary}
 \begin{proof}
 By Theorem \ref{thm-Serre}, the general extension is a vector bundle. Suppose $\mathcal{N} \to \cV$ is a destabilizing line subbundle. Then there must exist a nonzero map $\mathcal{N} \to I_Z \otimes \mathcal{L}$. Hence, $\mathcal{L}\otimes \mathcal{N}^{-1}$ is effective and nontrivial. In particular, $H^0(X, I_Z \otimes \mathcal{L} \otimes \mathcal{N}^{-1}) \not= 0$ and yet $0 < H \cdot c_1\left(\mathcal{L} \otimes\mathcal{N}^{-1}\right) < H \cdot c_1(\mathcal{L})$, contrary to our assumptions. 
 \end{proof}
 
 In general, it is hard to characterize loci of points on $X$ that satisfy the Cayley-Bacharach property for $\mathcal{L} \otimes \omega_X$. However, there are easy conditions that ensure that it holds.  For instance, we may choose $Z$ general and $|Z|$ large, or we may choose $\mathcal{L}$ such that $H^0(X, \mathcal{L} \otimes \omega_X)=0$.

\subsection{The moduli space is not necessarily irreducible}
Unlike in the case of curves, the stack of coherent sheaves with a fixed Chern character on a surface is in general reducible. Consequently, one cannot construct a coherent sheaf with the expected cohomology and deduce that the general stable sheaf will have the expected cohomology. Furthermore even the substack of semistable or stable sheaves can be reducible. 

\begin{example}[Mestrano \cite{Mestrano}]
Let $X$ be a very general surface of degree $d= 6$ in $\PP^3$. By the Noether-Lefschetz Theorem,  $\Pic(X)= \ZZ H$, where $H$ is the hyperplane class. Let $Z$ be a  zero-dimensional scheme of length  $11$ of one of the following types:
\begin{enumerate}
\item[ I)] $Z$ is contained in the intersection of $X$ with a twisted cubic curve $C$, or
\item[II)] $Z$ is a general set of points in a hyperplane section of $X$.
\end{enumerate}
 Consider an extension of the form 
$$0 \to \OO_X \to \cV \to I_Z(H) \to 0.$$ Since $\ext^1(I_Z(H), \OO_X) = h^1(X, I_Z(3H))=1$, up to scaling there exist unique nonsplit  extensions of this form.  Since $K_X = 2H$, to apply Theorem \ref{thm-Serre} we need to check that $Z$ satisfies the Cayley-Bacharach property with respect to $\OO_X(3H)$. By B\'ezout's Theorem a  hypersurface of degree $3$ that contains at least $10$ points of $C$ contains all of $C$. Hence, $Z$ satisfies the Cayley-Bacharach property in Case I. Similarly, in Case II,  any hypersurface of degree $3$ that contains 10  of the points contains the plane spanned by the points. Hence, $Z$ satisfies the Cayley-Bacharach property in Case II. By Theorem \ref{thm-Serre}, there are locally free sheaves $\cV$ in both cases. Observe that $\cV$ is necessarily $\mu_H$-stable since $\cV$ cannot admit a map from $\OO_X(kH)$ for $k \geq 1$. 

By a monodromy argument, one can check that the locus of bundles in Case I is irreducible of dimension 12, corresponding to a generically finite cover of the space of twisted cubics in $\PP^3$. The tangent space to the moduli space at $\cV$ is given by $\Ext^1(\cV, \cV) \cong H^1(X, \cV\otimes \cV^*)$. Since $\cV$ is rank 2, we have $\cV^* \cong \cV(-H)$ and consequently 
$$\cV \otimes \cV^* \cong (\cV \otimes \cV)(-H) \cong (\Sym^2 \cV \oplus \bigwedge^2 \cV)(-H) \cong \Sym^2(\cV)(-H) \oplus \OO_X.$$ Since $H^1(X, \OO_X)=0$, $\Ext^1(\cV, \cV) = H^1\left(X, \Sym^2(\cV)(-H)\right)$. Using the standard exact sequence for symmetric powers
$$0 \to \cV(-H) \to  \Sym^2(\cV)(-H) \to I_{2Z \subset X} (H) \to 0,$$ where $I_{2Z \subset X}$ is the symbolic square of the ideal of $Z$ in $X$, one shows that $\ext^1(\cV, \cV)=12$. Hence,  the bundles in Case I lie on a generically smooth component of the moduli space of dimension 12.

On the other hand, the locus of bundles in Case II has dimension at least 13. The choice of $Z$ depends on 14 parameters, 3 for the choice of a hyperplane $\Lambda$ and 11 for the choice of points in $\Lambda \cap X$. To go in the reverse direction, note that the bundles in Case II have $h^0(X, \cV) =2$, and from $\mu_H$-stability it follows that a choice of a nonzero global section determines a quotient $\cV/\OO_X\cong I_Z(H)$ such that $Z$ is in Case II from $h^0(X,\cV)=2$ and $h^1(X,\OO_X)=0$.  Therefore, the locus in the moduli space coming from the bundles in Case I has dimension at least 13. We conclude that the moduli space has at least 2 irreducible components.
\end{example}

By using a similar construction, one can produce examples of moduli spaces with arbitrarily many components parameterizing sheaves on surfaces in $\PP^3$.
\begin{theorem}\cite{CoskunHuizengaAbel}
For any integer $k$, there exists an integer $d_k$ such that if $d \geq d_k$, then a very general surface $X \subset \PP^3$ of degree $d$ has a moduli space of rank $2$ sheaves with $c_1 = H$ that has at least $k$ irreducible components. 
\end{theorem}

The moduli space does not even have to be connected. Okonek and Van de Ven \cite{Okonek} and Kotschick  \cite{Kotschick} found examples of disconnected moduli spaces on elliptic surfaces with high Picard rank. These examples use ample classes that are very close to the fiber class. One can even find disconnected examples on complete intersection surfaces of Picard rank 1.

\begin{example}
Let $d_1 \geq 5$ and $d_2 > d_1$ be two integers. Let $D_1$ be a smooth hypersurface in $\PP^4$ that contains lines.  Let $D_2$ be a very general hypersurface of degree $d_2$. Let $X$ be the complete intersection of $D_1$ and $D_2$, which by  Noether-Lefschetz Theory is a surface of Picard rank 1. Take a plane $\Lambda$ containing a line $\ell$ of $D_1$. Then $\Lambda \cap D_1 =  C \cup \ell$, where $C$ is the residual curve of degree $d_1 -1$. Let $Z = C \cap D_2$ be the zero dimensional subscheme of length $(d_1-1)d_2$. Then $Z$ satisfies the Cayley-Bacharach property for $\omega_X(H)$ and up to scalars there is a unique nonsplit extension of the form
$$0 \to \OO_X \to \cV \to I_Z(H) \to 0.$$  One can compute that $h^0(X, \cV) = 3$ and $h^1(X, \cV)=0$. Given a bundle $\cW$ in the same irreducible component as $\cV$, we have $h^0(X, \cW) \geq 3$. The cokernel of the evaluation map of  a section is $I_{Z'}(H)$, where $Z'$ has length $(d_1-1)d_2$. Furthermore, since $h^0(X, I_{Z'}(H))\geq 2$, $Z'$ lies on a plane $\Lambda'$.  By the Cayley-Bacharach Theorem, the residual $d_2$ points must be collinear and by B\'{e}zout's Theorem the line must be contained in $D_1$. With a slightly more careful analysis, one can deduce the following.
\end{example}

\begin{theorem}\cite[Theorem 3.4]{CoskunHuizengaKopperDisconnected}
Let ${\bf v}$ be the Chern character of the sheaf $\cV$ constructed in this example.
For every connected (respectively, irreducible) component of the Fano scheme of lines $F_1(D_1)$ on $D_1$, $M_{X,H}({\bf v})$ has a connected (respectively, irreducible) component of the same dimension.
\end{theorem}

By letting the degree of $D_1$ tend to infinity, we can find threefolds $D_1$ in $\PP^4$ that contain arbitrarily many  isolated lines. We thus deduce the following.

\begin{corollary}\cite[Corollary 1.1]{CoskunHuizengaKopperDisconnected}
For any integer $k$, there exists a smooth complete intersection surface $X \subset \PP^4$ with Picard rank 1 and a Chern character ${\bf v}$ on $X$ such that $M_{X,H}({\bf v})$ has at least $k$ connected components.
\end{corollary}
It is also interesting to note that the bundles $\cV$ do not deform to bundles on complete intersections where $D_1$ does not contain a line.

Even on rational surfaces there are disconnected moduli spaces.  Assuming the SHGH Conjecture \ref{conj-SHGH}, one can construct moduli spaces with arbitrarily many components for certain special polarizations (see \cite{CoskunHuizengaBlowupP2}).

In this subsection, we have discussed the connectedness and the irreducibility of the moduli spaces. In general, there are many interesting questions about the topology of $M_{X, H}({\bf v})$. A general  expectation due to Donaldson, Gieseker and Li is that the moduli spaces become better behaved as $\Delta$ becomes large. O'Grady's Theorem \ref{thm-OGrady} 
proves that the moduli space becomes irreducible for large $\Delta$. There are several conjectures concerning the stabilization of other  Betti numbers of the moduli spaces. 

\begin{conjecture}\cite[Conjecture 1.1]{CoskunWoolf}
Let $X$ be a smooth projective surface and let $H$ be an ample line bundle. Fix a rank $r > 0$ and a first Chern class $c$. Then the $i$th Betti number of $M_{X, H}(r, c, \Delta)$ stabilizes to a constant $b_{i, \Stab}(X)$, which is independent of $r$, $c$ and $H$, as $\Delta$ tends to $\infty$. 
\end{conjecture}
The conjecture is known when $X$ is a K3 or abelian surface and $r$ and $c$ are relatively prime. The conjecture is also known when $X$ is a rational surface, $M_{X,H}(r, c, \Delta)$ do not contain strictly semistable sheaves and $H$ is a polarization such that $H \cdot K_X < 0$ (see \cite{CoskunWoolf} for a discussion and references). The conjecture is wide open in general, especially for surfaces of general type.

\section{Weak Brill-Noether}\label{sec-WBN}
In this section, we will survey some results on the weak Brill-Noether Problem. Recall that  an irreducible component of a moduli space of sheaves satisfies {\em weak Brill-Noether} if the general sheaf in that component has at most one nonzero cohomology group. In particular, if $\chi({\bf v}) =0$ and a component of the moduli space $M({\bf v})$ satisfies weak Brill-Noether, then the general member of that component has no cohomology. The ideal situation occurs for  $\PP^2$. Let $L$ be the class of a line on $\PP^2$.

\begin{theorem}[G\"{o}ttsche-Hirschowitz \cite{GottscheHirschowitz}]\label{thm-GottscheHirschowitz}
A general stable bundle on $\PP^2$ has at most one nonzero cohomology group.
\end{theorem}
Consequently, the slope and the Euler characteristic determine the cohomology of the general stable bundle $\cV \in M_{\PP^2, L}({\bf v})$.
\begin{itemize}
\item If $\chi(\cV) <0$, then $h^1(\PP^2, \cV) = - \chi(\cV)$ and $h^0(\PP^2, \cV) = h^2(\PP^2, \cV) =0$. 
\item If $\chi(\cV) \geq  0$ and $\mu_L(\cV) \geq 0$, then $h^0(\PP^2, \cV) = \chi(\cV)$ and all other cohomology vanishes.
\item If $\chi(\cV) \geq 0$ and $\mu_L(\cV) < 0$, then $h^2(\PP^2, \cV) = \chi(\cV)$ and all other cohomology vanishes.
\item In particular, if $\chi(\cV) = 0$, then the cohomology of $\cV$ vanishes.
 \end{itemize}
Observe that if $\mu_L(\cV)< 0$, then $h^0(\PP^2, \cV) =0$ by stability. By Serre duality and stability, if $\mu_L(\cV) > -3$, then $h^2(\PP^2, \cV) =0$.

\begin{remark}
The sheaves $I_Z(-d)$ have nonvanishing $h^1$ and $h^2$ provided $|Z| \geq 1$ and $d \geq 3$. The general element of the corresponding moduli spaces are not locally free and Serre duality fails. The general stable sheaf of rank at least 2 on $\PP^2$ is locally free, hence weak Brill-Noether holds for $M_{\PP^2, L} ({\bf v})$ if the rank of ${\bf v}$ is at least 2.
\end{remark}

We will discuss two general techniques for proving weak Brill-Noether Theorems on surfaces. The first uses prioritary sheaves and the second uses Bridgeland stability conditions. We begin by explaining the first.

\subsection{Prioritary Sheaves} Let $D$ be an effective divisor on a smooth surface $X$. A torsion free sheaf $\cV$ is called {\em $D$-prioritary} if $\Ext^2(\cV, \cV(-D))=0$. We denote the stack of $D$-prioritary sheaves on $X$ with Chern character ${\bf v}$ by $\cP_{X,D}({\bf v})$. The stack $\cP_{X, D}({\bf v})$ is an open substack of the stack of coherent sheaves on $X$ with Chern character ${\bf v}$.

The condition $\Ext^2(\cV, \cV(-D))=0$ implies that the restriction map $$\Ext^1_X(\cV, \cV) \to \Ext^1_D(\cV|_D, \cV|_D)$$ is surjective. Consequently, the general first order deformation of the sheaf $\cV$ on $X$ gives a general first order deformation of $\cV|_D$ on $D$. For example, when $D$ is a smooth rational curve and $\cV$ is locally free on $D$, the restriction of a general deformation of $\cV$  to $D$ is balanced (i.e., the splitting type $\cV|_D = \oplus_{i=1}^r \OO_{\PP^1}(a_i)$ satisfies $|a_i - a_j| \leq 1$ for all $i,j$).

Let $H$ be an ample divisor on $X$.  Suppose that $(K_X+D) \cdot H <0$. Then a $\mu_H$-semistable sheaf $\cV$ is $D$-prioritary. By Serre duality,
$$\ext^2 (\cV, \cV(-D))= \hom(\cV, \cV(K_X+D)).$$ Since $(K_X + D) \cdot H <0$, stability implies that the latter group vanishes. Hence, $\cV$ is $D$-prioritary.

\begin{remark}
The condition $(K_X + D) \cdot H < 0$ implies that $K_X \cdot H < 0$. Hence, $mK_X$ cannot have any sections for $m > 0$. Therefore, this condition can only be satisfied for surfaces of Kodaira dimension $-\infty$. Furthermore, if $K_X + D$ is effective, there cannot be any $D$ prioritary sheaves on $X$. 
\end{remark}

A projective surface $X$  is {\em ruled} if it admits a morphism $X \to C$ onto a smooth curve where all the fibers are isomorphic to $\PP^1$, equivalently if $X$ is the projectivization of a vector bundle of rank 2 on $C$.  A surface is {\em birationally ruled} if it admits a morphism $X \to C$ onto a smooth curve where the general fiber is $\PP^1$. The following theorem of Walter makes prioritary sheaves a useful tool for studying moduli spaces of sheaves on birationally ruled surfaces. 

\begin{theorem}[Walter \cite{Walter}]
Let $X$ be a birationally ruled surface and let $F$ be the fiber class. Then the stack of $F$-prioritary sheaves $\mathcal{P}_F({\bf v})$ is irreducible whenever it is nonempty.
\end{theorem} 
This theorem generalizes an earlier theorem of Hirschowitz and Laszlo \cite{HirschowitzLaszlo} which assets that  $\mathcal{P}_{\PP^2, L} ({\bf v})$ is irreducible when nonempty. The main advantage of working with prioritary sheaves is that they are easier than stable sheaves to construct.  For example, $\OO_{\PP^2}(a) \oplus \OO_{\PP^2}(a+1)$ is $L$-prioritary on $\PP^2$ but not stable.

\begin{remark}
For rational surfaces, there are always polarizations $H$ for which $(K_X + F) \cdot H <0$. For minimal rational surfaces $K_X + F$ is anti-effective. Hence, $(K_X + F)\cdot H < 0$ for every ample class $H$. Every rational surface is obtained by blowing up a minimal rational surface at finitely many (possibly infinitely near) points. Let $Y$ be the blowup of $X$ at a point. If  $(K_X + F)\cdot H < 0$, then  $K_Y = \pi^* K_X + E$ and $$(K_Y + F) \cdot \pi^* H = (K_X + F) \cdot H <0.$$ The divisor $\pi^* H $ is not ample on $Y$, but nef.  Any small perturbation $H'$ of $\pi^* H$  by adding an ample divisor is ample. Hence, for an ample $H'$ sufficiently close to $\pi^* H$, we still have $(K_Y + F) \cdot H' < 0$.  
\end{remark}

The utility of the notion of $F$-prioritary sheaves toward the weak Brill-Noether problem follows from openness of (semi)stability and semi-continuity of cohomology.  Indeed, if $(K_X + F) \cdot H <0$ and there exist (semi)stable sheaves with Chern character ${\bf v}$, then to show that ${\bf v}$ satisfies weak Brill-Noether it suffices to show that there exists an $F$-prioritary sheaf with at most one nonzero cohomology group. By semi-continuity of cohomology and irreducibility of the moduli stack, it follows that the general stable sheaf has at most one nonzero cohomology group. Hence, the moduli space satisfies weak Brill-Noether. This strategy can be used to show weak Brill-Noether on certain rational or birationally ruled surfaces. 

\subsection{Elementary modifications} Given a torsion free sheaf $\mathcal{V}$, a point $p\in X$ and a surjection $\phi: \mathcal{V} \twoheadrightarrow \OO_p$, the kernel $\mathcal{V}'$ defined by the exact sequence 
\begin{equation}\label{eq-elementarymodification}
0 \to \mathcal{V}' \to \mathcal{V} \stackrel{\phi}{\to} \OO_p \to 0
\end{equation}
is called the {\em elementary modification} of $\mathcal{V}$ with respect to $\phi$. An elementary modification satisfies
$$r(\mathcal{V}') = r(\mathcal{V}), \quad \nu(\mathcal{V}') = \nu(\mathcal{V}), \quad \Delta(\mathcal{V}') = \Delta(\mathcal{V}) + \frac{1}{r}.$$
If $\mathcal{W}$ is a proper subsheaf of $\mathcal{V}'$ of smaller rank with $\mu_H(\mathcal{W}) \geq \mu_H(\mathcal{V}')$, then  $\mathcal{W}$ is also a proper subsheaf of $\mathcal{V}$ with $\mu_H(\mathcal{W}) \geq \mu_H(\mathcal{V})$. Hence, if $\mathcal{V}$ is a $\mu_H$-(semi)stable sheaf, then $\mathcal{V}'$ is also $\mu_H$-(semi)stable. If $\mathcal{V}$ is a $D$-prioritary sheaf, then a general elementary modification is also $D$-prioritary.

\begin{caution}
Elementary modifications do not necessarily  preserve Gieseker (semi)stability. For example, on $\PP^2$ the sheaf $\OO_{\PP^2} \oplus \OO_{\PP^2}$ is Gieseker semistable, however any elementary modification is isomorphic to $I_p \oplus \OO_{\PP^2}$, which is not Gieseker semistable.  
\end{caution}

Since $h^i(\OO_p)=0$ for $i>0$, we have that $H^2(X, \mathcal{V}') = H^2(X, \mathcal{V})$. The morphism $\phi$ corresponds to a choice of hyperplane $\Lambda$ in the fiber of $\mathcal{V}$ over $p$. If $h^0(X, \mathcal{V}) > 0$, we can choose $p$ and the hyperplane $\Lambda$ so that one of the sections at $p$ is not contained in $\Lambda$. In that case, the map $H^0(X, \mathcal{V}) \to H^0(X, \OO_p)$ induced by  $\phi$  is surjective. By the long exact sequence of cohomology associated to \eqref{eq-elementarymodification}, we conclude that 
$$h^0(X, \mathcal{V}') = h^0(X, \mathcal{V})-1, \quad h^1(X, \mathcal{V}') =  h^1(X, \mathcal{V}).$$

The integrality of the Euler characteristic and the Riemann-Roch Theorem
$$\chi(\mathcal{V})= r(\mathcal{V}) (P(\nu(\mathcal{V}))- \Delta(\mathcal{V}))$$
 imply that the difference of the discriminants of any two sheaves with the same rank and first Chern class is a multiple of $\frac{1}{r}$. Our discussion on elementary modifications yields the following theorem.
 
 \begin{theorem}\label{thm-reductiontosmalldiscriminant}
 Let $\mathcal{V}$ be a $\mu_H$-(semi)stable   (respectively, $D$-prioritary) sheaf of rank $r$, total slope $\nu$ and discriminant $\Delta_0$ on a smooth projective surface $X$ such that  $h^2(X, \mathcal{V})=0$ and $\mathcal{V}$ has at most one nonzero cohomology group. Then there exists a $\mu_H$-(semi)stable (respectively, $D$-prioritary) sheaf with rank $r$, total slope $\nu$ with at most one nonzero cohomology group for every $\Delta \geq \Delta_0$ for which $(r, \nu, \Delta)$ is an integral Chern character.
 \end{theorem}

\subsection{Weak Brill-Noether on $\PP^2$} We are now ready to sketch the proof of Theorem \ref{thm-GottscheHirschowitz}]  We will prove the following more general statement.
\begin{theorem}
Let $\cV$ be a general prioritary sheaf  on $\PP^2$ with $\Delta(\cV) \geq 0$ and rank at least 2. Then $\cV$ has at most one nonzero cohomology group. 
\end{theorem}
\begin{proof}
Since stable sheaves are prioritary and form a dense open substack of the stack of prioritary sheaves when nonempty, Theorem \ref{thm-GottscheHirschowitz} follows by the semicontinuity of cohomology. By the irreducibility of the stack of prioritary sheaves and semicontinuity of cohomology, it suffices to exhibit one prioritary sheaf with at most one nonzero cohomology group. When the rank is at least 2, the general prioritary sheaf is locally free. Hence, by Serre duality, we may assume that $\mu(\cV) \geq - \frac{3}{2}$.  Given a rank $r$ and a slope $\frac{c}{r} \geq - \frac{3}{2}$, we can find an $L$-prioritary sheaf $\cV= \OO_{\PP^2}(a)^{\oplus r-s} \oplus \OO_{\PP^2}(a+1)^{\oplus s}$ with $a \geq -2$ of rank $r$ and slope $\frac{c}{r}$. An easy computation shows that $\Delta(\cV) \leq 0$. The sheaf $\cV$ is nonspecial. Taking general elementary modifications, we obtain a prioritary sheaf that has at most one nonzero cohomology group for every integral Chern character.  This concludes the proof of the theorem.
\end{proof}

\subsubsection{Gaeta resolutions}  On $\PP^2$ one can write down a resolution for the general prioritary sheaf $\cV$ with $\Delta \geq 0$. By Riemann-Roch, there exists a unique integer $n$ such that $\chi(\cV(-n)) \geq 0$ but $\chi(\cV(-n-1)) < 0$. Set $n=a+2$. Then $\cV$ has a resolution of the form
 $$0 \to \OO_{\PP^2}(a)^{\alpha} \to \OO_{\PP^2}(a+1)^{\beta}\oplus \OO_{\PP^2}(a+2)^{\gamma} \to \cV \to 0, \ \mbox{or}$$
 $$0 \to \OO_{\PP^2}(a)^{\alpha} \oplus \ \OO_{\PP^2}(a+1)^{\beta}  \to \OO_{\PP^2}(a+2)^{\gamma} \to \cV \to 0.$$
 Here  $$\gamma = \chi(\cV(-a-2)), \quad \alpha = - \chi(\cV(-a-3)), \quad \beta = |\gamma - r(\cV) - \alpha |.$$ The resolution is of the first type if $\gamma - r(\cV) - \alpha \geq 0$ and otherwise of the second type. This generalizes Gaeta's resolution for the ideal sheaf of a general set of points on $\PP^2$ \cite{Gaeta}.
 
To see that such a resolution exists, given a Chern character ${\bf v}$ one can solve for $a, \alpha, \beta$ and $\gamma$ and write a general Gaeta-type resolution. One then checks that $\cV$ defined by such a resolution is prioritary and that the associated Kodaira-Spencer map is surjective  (see \cite[\S 3]{CoskunHuizengaWBN} or \cite[\S 4]{CoskunHuizengaBN} for more details). When $r(\cV) \geq 2$, by a Bertini-type theorem, a sheaf defined by such a resolution is locally free. Consequently, we also deduce that the general prioritary sheaf of rank at least 2 is locally free. The existence of the Gaeta resolution also implies that the moduli spaces $M_{\PP^2, L}({\bf v})$ are unirational.
 
 In general, we do not know a presentation for the general stable bundle on a surface. Having such presentations even for special families of surfaces would be very useful.

\subsubsection{Global generation} The weak Brill-Noether Theorem allows one to classify Chern characters of nonspecial stable bundles that are globally generated.

\begin{theorem}\cite[Corollary 5.3]{CoskunHuizengaBN} $($see also \cite{BGJ}$)$
The general member of $M_{\PP^2, L}({\bf v})$ is globally generated if and only if one of the following holds:
\begin{itemize}
\item ${\bf v} = r (1,0,0)$ and $\cV = \OO_{\PP^2}^{\oplus r}$
\item $\mu({\bf v}) > 0$ and $\chi({\bf v} (-1)) \geq 0$
\item $\mu({\bf v}) > 0$, $\chi({\bf v} (-1)) < 0$ and $\chi({\bf v}) \geq r+2$
\item $\mu({\bf v}) > 0$, $\chi({\bf v} (-1)) < 0$,  $\chi({\bf v}) \geq r+1$, and ${\bf v} = (r + 1) \ch(\OO_{\PP^2}) - \ch(\OO_{\PP^2}(-2))$.
\end{itemize}
\end{theorem}

\begin{proof}
If $\cV$ is globally generated, then its determinant is  globally generated, hence $\mu(\cV) \geq 0$.  If $\mu(\cV)=0$, then the Riemann-Roch Theorem implies that $\chi(\cV) \leq \rk(\cV)$ with equality if and only if $\Delta (\cV) =0$. If $\mu(\cV) = \Delta(\cV)=0$, then $\cV = \OO_{\PP^2}^{\oplus r}$.

If $\chi({\bf v}(-1)) \geq 0$, then the general sheaf in $\cP_{\PP^2, L}({\bf v})$ has a Gaeta resolution with $a \geq 1$. Then the general sheaf is clearly a quotient of a globally generated bundle.  If $\chi({\bf v}(-1)) < 0$ and $\chi({\bf v }) \geq \rk({\bf v}) + 2$, then the general sheaf in $\cP_{\PP^2, L}({\bf v})$ has a Gaeta resolution of the form 
$$0 \rightarrow \OO_{\PP^2}(-2)^k \oplus \OO_{\PP^2}(-1)^l \rightarrow \OO_{\PP^2}^m \rightarrow \cV \rightarrow 0, \  \mbox{or}$$ 
$$0 \rightarrow \OO_{\PP^2}(-2)^k  \rightarrow  \OO_{\PP^2}(-1)^l \oplus \OO_{\PP^2}^m \rightarrow  \cV \rightarrow 0.$$
In the first case, $\cV$ is the quotient of a globally generated vector bundle, hence globally generated. The most interesting case is the second one. By the assumption that $\chi(\cV) \geq \rk(\cV) + 2$, we have that $m \geq \rk(\cV) + 2$. Therefore, $k \geq l+2$. To show that $\cV$ is globally generated, it suffices to show that $H^1(\PP^2, \cV \otimes I_p) =0$ for every point $p \in \PP^2$. By the long exact sequence of cohomology, it suffices to show that the map 
$$\phi: H^1(\PP^2, I_p(-2))^k \rightarrow H^1(\PP^2, I_p(-1))^l$$ is surjective.
Consider the sequence $$0 \rightarrow M \rightarrow \OO_{\PP^2}(-2)^k  \rightarrow  \OO_{\PP^2}(-1)^l \rightarrow 0.$$ Since the map is general, it is surjective and $M$ is a vector bundle. Clearly $M$ does not have any cohomology. Tensoring the standard exact sequence $0 \rightarrow I_p \rightarrow \OO_{\PP^2} \rightarrow \OO_p \rightarrow 0$ with $M$, we see that $H^2(\PP^2, I_p \otimes M)=0$. Consequently, the map $\phi$ is surjective and $\cV$ is globally generated.
Finally, if $\chi(\cV) =\rk(\cV) +1$ and $\cV$ is globally generated, then there is a surjective map $\OO_{\PP^2}^{r+1} \rightarrow \cV$. The kernel of this map is a line bundle $\OO_{\PP^2}(-d)$. If $d=1$, then $\chi(\cV(-1)) =0$. If $d \geq 3$, then $\chi(\cV) < r$ and it is not possible for the general prioritary sheaf with Chern character ${\bf v}$ to be globally generated.  The only remaining possibility is for $d=2$. In that case, $\chi(\cV) = r+1$ and this is the Gaeta resolution of the general sheaf. This concludes the classification of globally generated Chern characters on $\PP^2$. 
 \end{proof}

\subsection{Weak Brill-Noether for Hirzebruch surfaces}
 Let $e \geq 0$ be an integer. Let  $\FF_e : = \PP(\OO_{\PP^1} \oplus \OO_{\PP^1}(e))$ denote the Hirzebruch surface. We refer the reader to \cite{Beauville, CoskunScroll} for  detailed information on the geometry of Hirzebruch surfaces. The surface $\FF_e$ has a section $E$ with self-intersection $-e$. Let $F$ denote the class of a fiber. Then   $\Pic(\FF_e) = \ZZ F \oplus \ZZ E$. The intersection numbers are $$E^2=-e, \quad E \cdot F = 1, \quad \mbox{and} \quad F^2 =0.$$ The effective cone of $\FF_e$ is spanned by $E$ and $F$. The nef cone of $\FF_e$ is spanned by $E+eF$ and $F$. The minimal rational surfaces are $\PP^2$ and $\FF_e$ for $e \geq 0$, $e\not= 1$.  The existence of curves with negative self-intersection provides an obstruction for  the weak Brill-Noether property on Hirzebruch surfaces $\FF_e$ with $e \geq 1$.

\subsubsection{Negative curves on a surface and the cohomology of vector bundles} Suppose that a surface $X$ contains a smooth curve $C$ with $C^2 < 0$. Then by the Riemann-Roch Theorem, $$\chi(\OO_X(mC))= \chi(\OO_X) + \frac{m^2}{2} C^2 - \frac{m}{2} K \cdot C.$$ If $m \gg 0$, $\chi(\OO_X(mC)) < 0$. We conclude that both $H^0(X, \OO_X(mC))$ and $H^1(X, \OO_X(mC))$ are nonzero. In fact, if $X$ is a surface such that  $$q=h^1(X, \OO_X)=0 \quad \mbox{and} \quad p_g = h^2(X, \OO_X)=0,$$ then  $\OO_X(mC)$ is  an effective line bundle with nonvanishing $h^1$ provided that either
\begin{enumerate}
\item  $m \geq 1$ and either  $C^2 < -1$ or the genus of $C$ is at least $1$, or  
\item $m \geq 2$ and $C$ is an exceptional curve. 
\end{enumerate}
A similar phenomenon persists for higher rank vector bundles.

\begin{proposition}\label{prop-wBNfailsnegative}
Let $X$ be a smooth projective surface that contains an irreducible curve $C$ with $C^2 < 0$. Then for every rank $r$ there exist infinitely many Chern characters ${\bf v}_i$ with rank $r$ such that the moduli space $M_{X, H}({\bf v}_i)$ is nonempty and does not satisfy weak Brill-Noether. 
\end{proposition}

\begin{proof}
Let $M_{X, H}({\bf v})$ be an irreducible moduli space of rank $r$ sheaves where the general  sheaf is $\mu_H$-stable and has sections. Such moduli spaces exist by O'Grady's Theorem \ref{thm-OGrady} and Serre vanishing. Let $\cV$ be a general sheaf in  $M_{X, H}({\bf v})$. Let ${\bf v}_i$ be the Chern character of $\cV(iC)$. Since $H^0(X, \cV) \subset H^0(X, \cV(iC))$, the sheaves $\cV(iC)$ have global sections. On the other hand, by Riemann-Roch if $i \gg 0$, then $\chi({\bf v}_i) < 0$, hence, they must also have $h^1$ and cannot satisfy weak Brill-Noether. 
\end{proof}

Unlike $\PP^2$, by Proposition \ref{prop-wBNfailsnegative}, we cannot expect weak Brill-Noether to always hold for moduli spaces on $\FF_e$ with $e \geq 1$. The next theorem shows that the existence of negative curves is the only obstruction for weak Brill-Noether on $\FF_e$.

\begin{theorem}\cite[Theorem 3.1]{CoskunHuizengaBN}\label{thm-HirzebruchWBN}
Let ${\bf v}$ be a Chern character with positive rank $r({\bf v})$ and $\Delta({\bf v}) \geq 0$. Then the stack of prioritary sheaves $\mathcal{P}_F({\bf v})$ is nonempty and irreducible. Let $V$ be a general sheaf in $\mathcal{P}_F({\bf v})$.
\begin{enumerate}
\item If $\nu({\bf v}) \cdot F \geq -1$, then $h^2(\FF_e, \mathcal{V}) =0$.
\item If $\nu({\bf v}) \cdot F \leq -1$, then $h^0(\FF_e, \mathcal{V}) =0$.
\item If $\nu({\bf v}) \cdot F = -1$, then $h^1(\FF_e, \mathcal{V}) = - \chi(\mathcal{V})$ and all other cohomology vanishes.
\end{enumerate}
Assume that $\nu({\bf v}) \cdot F > -1$.
\begin{enumerate}
\item[(4)] If $\nu({\bf v}) \cdot E \geq -1$, then $\mathcal{V}$ has at most one nonzero cohomology group. If $\chi(\mathcal{V}) \geq 0$, then $h^0(\FF_e, \mathcal{V}) =\chi(\mathcal{V})$, and if $\chi(\mathcal{V}) < 0$, then $h^1(\FF_e, \mathcal{V}) =-\chi(\mathcal{V})$.
\item[(5)] If $\nu({\bf v}) \cdot E < -1$, then $H^0(\FF_e, \mathcal{V}) \cong H^0(\FF_e, \mathcal{V}(-E))$ and the Betti numbers of $\mathcal{V}$ are inductively determined by (2) and (4).
\end{enumerate}
\end{theorem}
Observe that if $\nu({\bf v}) \cdot F < -1$ and $\rk({\bf v}) \geq 2$, then the general prioritary sheaf is locally free and Serre duality determines the cohomology.

\begin{proof}[Sketch of proof]
The proof is similar to the proof of Theorem \ref{thm-GottscheHirschowitz}. A vector bundle $\cV$ of the form 
$$ \OO_{\FF_e}(-E-(e+1)F)^{\oplus a} \oplus \OO_{\FF_e}(-F)^{\oplus b} \oplus \OO_{\FF_e}^{\oplus c} \ \mbox{or} \  \OO_{\FF_e}(-E-(e+1)F)^{\oplus a} \oplus \OO_{\FF_e}(-E-eF)^{\oplus b} \oplus \OO_{\FF_e}^{\oplus c}$$
is both $E$-prioritary and $F$-prioritary and has $\Delta(\cV) \leq 0$. Moreover, $\cV$ has no higher cohomology. If ${\bf v}$ satisfies $\nu({\bf v}) \cdot F \geq -1$ and $\nu({\bf v}) \cdot E \geq -1$, we can find a prioritary bundle of Chern character ${\bf v}$ by tensoring ${\bf v}$ by a nef line bundle $\mathcal{N}$ on $\FF_e$ since tensoring with a nef line bundle preserves prioritariness and the discriminant.   Furthermore, $\cV\otimes \mathcal{N}$ has no higher cohomology on $\FF_e$. We conclude that weak Brill-Noether holds for slopes satisfying $\nu({\bf v}) \cdot F \geq -1$ and $\nu({\bf v}) \cdot E \geq -1$. For slopes in the range $\nu({\bf v}) \cdot F > -1$ and $\nu({\bf v}) \cdot E <-1$, one can use the exact sequence
$$0 \to \cV(-E) \to \cV \to \cV|_E \to 0$$ and the fact that $\cV|_E$ is balanced on  $E$ to inductively compute the cohomology.
\end{proof}

As in the case of $\PP^2$, the weak Brill-Noether Theorem allows one to  find a Gaeta-type resolution for the general prioritary sheaf on Hirzebruch surfaces. 
\begin{theorem}[\cite{CoskunHuizengaBN}, Theorem 4.1]\label{thm-Gaeta}
Let ${\bf v}$ be an integral Chern character on $\FF_e$ of positive rank and assume that 
$$\Delta({\bf v}) \geq \frac{1}{4} \ \mbox{if} \ e=0, \quad \Delta({\bf v}) \geq \frac{1}{8} \ \mbox{if} \ e=1, \quad  \Delta({\bf v}) \geq 0 \ \mbox{if} \ e \geq 2.$$
Then the general sheaf $\cV \in \cP_{\FF_e, F}({\bf v})$ admits a Gaeta-type resolution 
\begin{equation}\label{eq-Gaeta}
0 \rightarrow L(-E-(e+1)F)^a \rightarrow L(-E-eF)^b \oplus L(-F)^c \oplus L^d \rightarrow \cV \rightarrow 0,
\end{equation}
for some line bundle $L$ and nonnegative integers $a, b, c, d$. 
\end{theorem}
Theorems \ref{thm-HirzebruchWBN} and \ref{thm-Gaeta} allow one  to classify moduli spaces where the general bundle is a nonspecial, globally generated bundle.

\begin{theorem}\cite[Theorem 5.1]{CoskunHuizengaBN}
Suppose $e \geq 1$. Let ${\bf v}$ be a Chern character on $\FF_e$ such that $r({\bf v}) \geq 2$, $\Delta({\bf v}) \geq 0$ and $\nu({\bf v})$ is nef. Then the general member of $\mathcal{P}_F({\bf v})$ is globally generated if and only if one of the following holds.
\begin{enumerate}
\item We have $\nu({\bf v}) \cdot F =0$ and $${\bf v}= (r({\bf v}) - m) \ch(\OO_{\FF_e}(aF)) + m  \ch(\OO_{\FF_e}((a+1)F))$$ for integers $a, m \geq 0$.
\item We have $\nu({\bf v}) \cdot F >0$ and $\chi({\bf v}(-F)) \geq 0$.
\item We have $\nu({\bf v}) \cdot F >0$ and $\chi({\bf v}(-F)) < 0$ and $\chi({\bf v}) \geq r({\bf v}) + 2$.
\item We have $e=1$,  $\nu({\bf v}) \cdot F >0$,  $\chi({\bf v}(-F)) < 0$, $\chi({\bf v}) \geq r({\bf v}) + 1$ and $${\bf v} = (r({\bf v} + 1) \ch\left( \OO_{\FF_1}\right) - \ch \left(\OO_{\FF_1}(-2E-2F)\right).$$
\end{enumerate}
\end{theorem}

\begin{remark}
For $\FF_0= \PP^1 \times \PP^1$,  the statement needs to be slightly modified (see \cite[Theorem 5.2]{CoskunHuizengaBN}). If $F_1$ and $F_2$ are the two rulings on $\PP^1 \times \PP^1$, then the general member of $\mathcal{P}_F({\bf v})$ is globally generated if and only if one of the following holds.
\begin{enumerate}
\item We have $\nu({\bf v}) \cdot F_i =0$ for some $1 \leq i \leq 2$ and $${\bf v}= (r({\bf v}) - m) \ch(\OO_{\FF_e}(aF_i)) + m  \ch(\OO_{\FF_e}((a+1)F_i))$$ for integers $a, m \geq 0$.
\item We have $\nu({\bf v}) \cdot F_i >0$ for $1 \leq i \leq 2$ and $\chi({\bf v}(-F_j)) \geq 0$ for some $1 \leq j \leq 2$.
\item We have $\nu({\bf v}) \cdot F_i >0$, $\chi({\bf v}(-F_i)) < 0$ for $1 \leq i \leq 2$ and $\chi({\bf v}) \geq r({\bf v}) + 2$.
\end{enumerate}
\end{remark}

\subsubsection{Applications to ample bundles} On curves stable vector bundles are ample if and only if their slope is positive. On higher dimensional varieties, it is easy to see that there can be no numerical characterization of ample vector bundles of higher rank. Ampleness is an open condition and we can ask when a general member of an irreducible component of a moduli space is ample.

\begin{problem}
Let $X$ be a smooth, projective variety. Classify Chern characters ${\bf v}$ on $X$ for which there exists an ample vector bundle of Chern character ${\bf v}$.
\end{problem}
This problem is wide open even for surfaces. Recently, Huizenga and Kopper \cite{HuizengaKopper} have made significant progress towards the solution of the problem for $\PP^2$ and Hirzebruch  surfaces. If $\cV$ is an ample bundle on $\PP^2$ or $\FF_e$, then the restrictions $\cV|_L$, $\cV|_F$ and $\cV|_E$ are ample, hence split as a direct sum of line bundles on $\PP^1$ of positive degree. Huizenga and Kopper prove that on minimal rational surfaces, asymptotically this is the only obstruction to ampleness of the general stable bundle.

\begin{theorem}\cite[Theorem 4.1]{HuizengaKopper}
Let $X= \PP^2$ or $\FF_e$. Let ${\bf v}$ be a stable Chern character on $X$ such that 
\begin{enumerate}
\item if $X= \PP^2$, $\nu({\bf v}) \cdot L > 1$,
\item if $X= \PP^1 \times \PP^1$, $\nu({\bf v}) \cdot F > 1$ and $\nu({\bf v}) \cdot E > 1$,
\item if $X= \FF_e$ with $e\geq 1$, $\nu({\bf v}) \cdot F > 1$ and $\nu({\bf v}) \cdot E \geq 1$.
\end{enumerate}
Then for $n \gg 0$, the general bundle $\cV \in M_{X, H} (n {\bf v})$ is ample.
\end{theorem}

Huizenga and Kopper also classify globally generated ample bundles on minimal rational surfaces.

\begin{theorem}\cite[Theorem 5.1]{HuizengaKopper}
Let $X= \PP^2$ or $\FF_e$. Let ${\bf v}$ be a Chern character such that the moduli space $M_{X, H} ({\bf v})$ is nonempty and the general sheaf is a globally generated vector bundle with no higher cohomology.
\begin{enumerate}
\item if $X= \PP^2$, $\nu({\bf v}) \cdot L > 1+ \frac{1}{r}$,
\item if $X= \PP^1 \times \PP^1$, $\nu({\bf v}) \cdot F > 1$ and $\nu({\bf v}) \cdot E > 1$,
\item if $X= \FF_e$ with $e\geq 1$, $\nu({\bf v}) \cdot F > 1$ and $\nu({\bf v}) \cdot E \geq 1$.
\end{enumerate}
Then the general sheaf $\cV \in M_{X, H}({\bf v})$ is ample.
\end{theorem}

\subsection{Weak Brill-Noether for more general rational surfaces} The ideas and techniques used to study the weak Brill-Noether Problem on $\PP^2$ and Hirzebruch surfaces can be extended to non-minimal rational surfaces. Here we will discuss general blowups of $\PP^2$ and refer the reader to \cite{CoskunHuizengaWBN, CoskunHuizengaAbel} for analogous statements for general blowups of Hirzebruch surfaces. 

Let $X_m$ denote the blowup of $\PP^2$ at $m$ general points $p_1, \dots, p_m$. Let $E_i$ denote the exceptional divisor lying over $p_i$. Let $H$ denote the pullback of the class of a line. Then $$\Pic(X_m) = \ZZ H \oplus \bigoplus_{i=1}^m E_i$$ and the intersection pairing satisfies $$H^2 =1, \quad H \cdot E_i = 0, \quad E_i \cdot E_j = - \delta_{i,j},$$ where  $\delta_{i,j}$ is the Kronecker delta function. When $m \leq 8$, then the corresponding surface is a del Pezzo surface. Del Pezzo surfaces and  $\PP^1 \times \PP^1$ are the Fano surfaces, that is the surfaces that have an ample anti-canonical bundle. We refer the reader to \cite{Beauville, CoskundelPezzo, Hartshorne} for more details on the geometry of del Pezzo surfaces.

Proposition \ref{prop-wBNfailsnegative} shows that we cannot expect weak Brill-Noether to always hold. However, the following higher rank generalization of the SHGH might hold.

\begin{conjecture}\cite[Conjecture 1.7]{CoskunHuizengaWBN}
Let $X$ be a blowup of $\PP^2$ at $m$ very general points. Let $H$ be an ample class on $X$ such that $H \cdot K_X < 0$. Let ${\bf v}$ be a Chern character such that $\nu({\bf v})$ is nef.  Then the general stable sheaf in $M_{X, H} ({\bf v})$ has at most one nonzero cohomology group.
\end{conjecture}

\begin{remark}
One can make a bolder conjecture by weakening the condition that   $\nu({\bf v})$  is nef to $C \cdot \nu({\bf v}) \geq -1 $ for every $(-1)$-curve on $X$.  This would mimic the SHGH Conjecture more closely. Even the weaker conjecture is open when $m \geq 5$. 
\end{remark}

When $\nu({\bf v})$ is not too close to the boundary of the nef cone, one can find prioritary direct sums of line bundles with no higher cohomology. More precisely, let ${\bf v}$ be a Chern character of rank $r$ and let the total slope  be $\nu({\bf v}) = \delta H - \sum_{i=1}^k \alpha_i E_i$. Then we have that $$\delta = d + \frac{q}{r}, \quad \alpha_i = a_i + \frac{q_i}{r}$$ for some integers $d, q, a_i$ and $q_i$ with $0 \leq q < r$ and $0 \leq q_i < r$. Set $$\gamma({\bf v}) = \frac{q^2}{2r^2} - \frac{q}{2r} + \sum_{i=1}^k \left(\frac{q_i}{2r}- \frac{q_i^2}{2r^2} \right).$$

\begin{theorem}\cite[Theorem 4.12]{CoskunHuizengaAbel}
Let $X$ be the blowup of $\PP^2$ at $m$ distinct points. Let ${\bf v}$ be a positive rank Chern character on $X$ with total slope 
$$\nu({\bf v}) = \delta H - \alpha_1 E_1 \cdots - \alpha_k E_k$$ with $\delta \geq 0$ and $\alpha_i \geq 0$. Suppose that the line bundle $$ \lfloor \delta \rfloor H -\lceil  \alpha_1 \rceil  E_1 \cdots - \lceil \alpha_k \rceil E_k$$ does not have higher cohomology.  Assume that $\Delta({\bf v}) \geq \gamma({\bf v})$.
Then the stack of prioritary sheaves $\cP_{X, H-E_1}({\bf v})$ is nonempty and the general sheaf in $\cP_{X, H-E_1}({\bf v})$ has at most one nonzero cohomology group.
\end{theorem}

Levine and Zhang in \cite{LevineZhang} have studied the case of del Pezzo surfaces in greater detail.  Let $W_m$ denote the Weyl group acting on $\Pic(X_m)$.   Let $\cP_{H^W}({\bf v})$ denote the stack parameterizing torsion free sheaves $\cV$ such that $\Ext^2(\cV, \cV(-\sigma(H)))=0$ for every $\sigma \in W_m$.
\begin{theorem}\cite[Theorem 1.2]{LevineZhang}
Let $X_m$ be a del Pezzo surface with $m \leq 5$. Let ${\bf v}$ be a Chern character such that $H \cdot \nu({\bf v}) \geq -2$ and $\cP_{H^W}({\bf v})\not= \emptyset$.
\begin{enumerate}
\item If $C \cdot \nu({\bf v}) \geq -1$ for all $(-1)$-curves $C$ on $X_m$, then ${\bf v}$ is nonspecial.
\item If there exists $\sigma \in W_m$ such that $\nu({\bf v}) \cdot \sigma(H) \leq -1$ or $\nu({\bf v}) \cdot \sigma(H-E_i) \leq -1$ for some $i$, then ${\bf v}$ is nonspecial.
\item Let $\nu({\bf v}) \cdot \sigma(H) >-1$ and $\nu({\bf v}) \cdot \sigma(H-E_i)> -1$ for all $i$ and all $\sigma \in W_m$. Suppose $C$ is a $(-1)$-curve such that $C \cdot \nu({\bf v}) < -1$ and let $\pi$ denote the map contracting $C$. Then ${\bf v}$ is nonspecial if and only if $\pi_*(\cV)$ is nonspecial for the general $\cV$ and $\chi(\pi_*(\cV))\leq 0$.
\end{enumerate}
\end{theorem}

As an application, when $m \leq 6$, Levine and Zhang classify $-K_{X_m}$-stable Chern characters in terms of a Dr\'{e}zet-Le Potier type condition. More generally, by studying blowups of $\PP^2$ along collinear points and taking deformations, Zhao has studied the weak Brill-Noether Problem on general blowups of $\PP^2$ \cite[Theorems 4.7 and 7.4]{Zhao} and has shown the nonemptiness of moduli spaces for certain Chern characters and polarizations \cite[Theorems 6.14 and 7.3]{Zhao}.

\subsection{Weak Brill-Noether for $K$-trivial surfaces}

There has been recent progress for computing the cohomology of the general sheaf on $K$-trivial surfaces using Bridgeland stability conditions. In this subsection, we will introduce Bridgeland stability conditions and briefly explain the strategy for using them to prove weak Brill-Noether Theorems, concentrating on the case of K3 surfaces.

\subsubsection{Bridgeland stability} Let $\cD^b(X)$ denote the bounded derived category of coherent sheaves on $X$. Let $K(\cD^b(X))$ denote the $K$-group of $\cD^b(X)$ and let $Z: K(\cD^b(X)) \to \CC$ be a group homomorphism called {\em the central charge} which factors through the Chern character. 

\begin{definition}
A {\em Bridgeland stability condition}  on $\cD^b(X)$  is a pair $\sigma= (\mathcal{A}, Z)$, where $\mathcal{A}$ is an abelian category which is the heart of a bounded $t$-structure on $\cD^b(X)$, and $Z$ is a central charge satisfying the following properties:
\begin{enumerate}
\item Positivity:  If $0 \not= E \in \mathcal{A}$, then $Z(E)=re^{i\theta}$ with $r>0$ and $0< \theta \leq \pi$.  Using $Z$, we  can define  $\mu_\sigma(E) = - \frac{\Re(Z(E))}{\Im(Z(E))}$. An object $E \in \mathcal{A}$ is {\em $\sigma$-(semi)stable} if for every subobject $F$ in $\mathcal{A}$, we have $\mu_{\sigma}(F) \leqor \mu_{\sigma}(E)$.
\item Harder-Narasimhan property: Every object $E \in \mathcal{A}$ has a finite Harder-Narasimhan filtration with $\sigma$-semistable quotients of strictly decreasing slope.  
\item Support condition: For a fixed norm $|\cdot |$ on $H^*_{\mbox{alg}}(X, \RR)$, there exists a constant $C > 0$ such that for all $\sigma$-semistable $E \in \mathcal{A}$ we have $|Z(E)|\geq C |\ch(E)|$.
\end{enumerate}
\end{definition}
The set $\Stab(X)$ of Bridgeland stability conditions on $X$ has the structure of a complex manifold \cite[Corollary 1.3]{Bridgeland}. Here we will be interested in very special stability conditions on surfaces constructed by Bridgeland \cite{BridgelandK3} and Arcara and Bertram \cite{ArcaraBertram}.

Given an ample divisor $H$ on a surface $X$,  define two subcategories of the category of coherent sheaves $\Coh(X)$ by
\begin{align*}\cT_{s} &= \{ E \in \Coh(X) : \mu_{H}(G) > s H^2\ \mbox{for every quotient} \ G \ \mbox{of} \ E\} \\
\cF_{s} &=\{ E \in \Coh(X) :  \mu_{H}(F) \leq s H^2 \ \mbox{for every subsheaf} \ F \ \mbox{of} \ E \}.\end{align*}
The pair $(\cT_{s}, \cF_{s})$ forms a torsion pair in $\Coh(X)$, that is $\Hom(T, F) = 0$ if $T \in \cT_s$ and $F \in \cF_s$ and every coherent sheaf $\cV$ fits in a unique exact sequence
$$0 \to T \to \cV \to F \to 0,$$ where $T \in \cT_s$ and $F \in \cF_s$. 

Given a torsion pair in the heart of a $t$-structure, one obtains the heart of a new $t$-structure by tilting. Explicitly, tilting $\Coh(X)$ with respect to $(\cT_s, \cF_s)$, we obtain the category $\cA_s$ defined by 
$$\cA_{s} = \{ E^{\bullet} \in \cD^b(X) : H^{-1}(E^{\bullet}) \in \cF_{s}, H^0(E^{\bullet}) \in \cT_{s}, H^i(E^{\bullet})=0 \ \mbox{for} \ i\not= -1,0\}.$$
 Define 
$$Z_{s,t} (E)= -\int_X e^{-(s+it)H} \ch(E).$$
Then the pair $\sigma_{s,t} = (\cA_{s}, Z_{s,t})$ is a Bridgeland stability condition for $s,t \in \RR, t >0$ \cite{ArcaraBertram}. Given a Chern character ${\bf v}$ and a stability condition $\sigma= (\cA_s, Z_{s,t})$, let $M_{\sigma}({\bf v})$ be the moduli space of $\sigma$-semistable objects in $\cA_s$ with Chern character ${\bf v}$.

 \subsubsection{Walls} Given a Chern character ${\bf v}$, there is a locally finite wall-and-chamber decomposition of the upper $(s,t)$-half-plane such that within each chamber the $\sigma$-(semi)stable objects remain constant. If a semistable object $E$ gets destabilized by a subobject $F$, then the $\sigma$-slopes of $E$ and $F$ must become equal for some stability conditions. Expressing the equality $\mu_{\sigma}(E) = \mu_{\sigma}(F)$ as a function of $s$ and $t$, one obtains that the walls are either vertical lines or nested, disjoint semi-circles.
 Moreover, when $t$ is sufficiently large (the large volume limit), the Bridgeland moduli space is isomorphic to $M_{X, H}({\bf v})$.

\subsubsection{Weak Brill-Noether for K3 surfaces} Let $X$ be a K3 surface. Stable sheaves on K3 surfaces have been classified by Mukai \cite{Mukai2}, O'Grady \cite{OGrady2} and Yoshioka \cite{Yoshioka99}. It is customary to phrase the classification in terms of the Mukai vector instead of the Chern character. Let 
$${\bf v}(E) := \ch(E) \sqrt{\mbox{td}(X)} = (r(E), c_1(E), r(E) + \ch_2(E)),$$ where $\mbox{td}(X)$ is the Todd class of $X$. There is a pairing between any two Mukai vectors given by 
$$\langle {\bf v}, {\bf v}' \rangle = \langle (r,c,a), (r', c', a') \rangle := c \cdot c' - ra' - r' a = - \chi({\bf v}, {\bf v}'), $$ where $c \cdot c'$ is the intersection product on $H^2(X, \ZZ)$. A class is called {\em spherical} if ${\bf v}^2 =-2$ and \emph{isotropic} if ${\bf v}^2 =0$. A Mukai vector $\bf v$ is \emph{primitive} if it is not divisible in $H^*_{\textrm{alg}} (X,\ZZ)$. We say a primitive Mukai vector ${\bf v} = (r,c,a)$ is \emph{positive} if ${\bf v}^2 \geq -2$ and either 
\begin{enumerate}
\item $r>0$; or
\item $r = 0$, $c$ is effective, and $a\neq 0$; or
\item $r=c=0$ and $a> 0$.
\end{enumerate}

\begin{theorem}\label{Thm:ClassicalModSp}
Let $X$ be a K3 surface and let ${\bf v}=m{\bf v}_0$ be a Mukai vector, where ${\bf v}_0$ is  a primitive positive Mukai vector and $m>0$.  Then $M_{X,H}({\bf v})$ is non-empty for any ample divisor $H$.  If $H$ is \emph{generic} with respect to ${\bf v}$, then:
\begin{enumerate}
    \item The  moduli space $M_{X,H}({\bf v})$ is non-empty if and only if ${\bf v}_0^2\geq-2$.
    \item If $m=1$ or ${\bf v}_0^2 >0$, then  $\dim M_{X,H}({\bf v})={\bf v}^2+2$.
    \item When ${\bf v}_0^2=-2$, then  $M_{X,H}({\bf v})$ is a single point parameterizing the direct sum of $m$ copies of a spherical bundle. When ${\bf v}_0^2=0$, then $\dim M_{X,H}({\bf v}) = 2m$.
    \item When ${\bf v}_0^2>0$, $M_{X,H} ({\bf v})$ is a normal irreducible projective variety with $\QQ$-factorial singularities.  
\end{enumerate}
\end{theorem}

Now suppose that  $\Pic(X) = \ZZ H$ with $H^2 = 2n$.  If $\cV$ is a stable sheaf with $\ch_1(\cV) = dH$ with $d >0$, then $H^2(X, \cV) =0$ by stability. Hence, we need to compute $H^0(X, \cV)$ and $H^1(X, \cV)$. The fact that $H^2(X, \OO_X) = \CC$ allows one to construct many counterexamples to weak Brill-Noether on K3 surfaces.

\begin{example}\cite[Example 1.3]{CoskunNuerYoshioka}
The linear system $|H|$ defines a morphism $f: X \to \PP^{n+1}$. The sheaf  $\cV = f^* T_{\PP^{n+1}}$ is the unique stable sheaf in its moduli space. The pullback of the Euler sequence by $f$ $$0 \to \OO_X \to \OO_X(H)^{n+2} \to \cV \to 0$$ implies that $h^0(X, \cV) = n^2 + 4n +3$ and $h^1(X, \cV) =1$.
\end{example}

\begin{example}\cite[Example 6.3]{CoskunNuerYoshioka}
Let $X$ be a double cover of $\PP^2$ branched along a very general sextic curve. Then the pullback of an exceptional bundle on $\PP^2$ is a spherical stable bundle on $X$,  which is the unique member of its moduli space.  For example, let $f_k$ be the $k$th Fibonacci number with $f_1 = f_2 =1$. For $k\geq 2$, the pullback of the corresponding exceptional bundle from $\PP^2$ has resolution
$$0 \to \OO_X^{\oplus f_{2k-2}} \to \OO_X(H)^{\oplus f_{2k}}\to \cV_k \to 0.$$ The bundle $\cV_k$ has rank $f_{2k-1}$ and $$h^0(X, \cV_k) = 2f_{2k} + f_{2k-1} \quad \mbox{and} \quad h^1(X, \cV_k) = f_{2k-2}.$$ When $k=2$, one recovers the $n=1$ case of the previous example.
\end{example}

\begin{example}\cite[Theorem 10.1]{CoskunNuerYoshioka}
Let $C$ be a smooth member in $|H|$. Consider extensions of the form 
$$0 \to \OO_X(H)^r \to \cV \to \OO_C(L) \to 0,$$ where $L$ is a general line bundle on $C$ with Euler characteristic $2n-r$.  Then $\v(\cV)=(r, (r+1)H, n(r+2))$, and from the exact sequence we conclude that $$h^0(X, \cV) = r(n+2) + \max(0, 2n-r) \quad \mbox{and} \quad h^1(X, \cV) = \max(0, r-2n).$$ As the general element in $M_{X,H}({\bf v})$ is given by such an extension, we get counterexamples to the weak Brill-Noether property when $r > 2n$.
\end{example}

These examples should convince the reader that the structure of the set of counterexamples to the weak Brill-Noether property on K3 surfaces is quite complicated and depends on arithmetic properties of the Mukai lattice.  We now explain the strategy to use Bridgeland stability conditions to prove weak Brill-Noether. 

\subsubsection{The strategy} Let $\Delta \subset X \times X$ be the diagonal and let $\pi_1$ and $\pi_2$ denote the two projections from $X\times X$ to $X$. Let  $\Phi_{X \rightarrow X}^{I_{\Delta}}: \cD^b(X) \to \cD^b(X)$ be the Fourier-Mukai transform acting by $\Phi_{X \rightarrow X}^{I_{\Delta}}(E) = \pi_{2*}(\pi_1^*(E)\otimes I_{\Delta})$. For each $E \in \Coh(X)$, tensoring the exact sequence $$ 0 \to I_{\Delta} \to \OO_{X\times X} \to \OO_{\Delta} \to 0$$ by $\pi_1^* E$ and pushing forward by $\pi_2$ gives the exact triangle
\begin{equation*}\label{eqn:FundamentalExactTriangle}
\Phi_{X\to X}^{I_\Delta}(E)\to R\Gamma(X,E)\otimes\OO_X\to E   
\end{equation*}
which induces the long exact sequence
$$ 0 \to  \mathcal{H}^0(\Phi_{X \to X}^{I_\Delta}(E)) \to  H^0(X,E) \otimes\OO_X \stackrel{f}{\to} E \to \mathcal{H}^1(\Phi_{X \to X}^{I_\Delta}(E)) \to  $$ $$H^1(X,E) \otimes \OO_X
\to  0\to  \mathcal{H}^2(\Phi_{X \to X}^{I_\Delta}(E)) \to H^2(X,E) \otimes \OO_X
\to 0.$$
Here $f$ is the evaluation morphism. Let $F:=  \Phi_{X \to X}^{I_\Delta}(E))^{\vee}$ be the derived dual. Then  $$\mathcal{H}^i(\Phi_{X \to X}^{I_\Delta}(E)))=\mathcal{H}^{i}(F^\vee)=\mathcal{E}xt^{i}(F,\OO_X).$$

\begin{lemma}\cite[Lemma 3.1]{CoskunNuerYoshioka}\label{lem-CNYmain}
Let $E$ be a coherent sheaf with no zero-dimensional torsion and set $F:=\Phi_{X \to X}^{I_\Delta}(E)^{\vee}$.  Then 
\begin{enumerate}
\item
$F$ is a coherent sheaf if and only if $E$ is nonspecial and generically globally generated.
\item
$F$ is a torsion-free sheaf if and only if $E$ is nonspecial and fails to be globally generated in at most finitely many points.
\item
$F$ is a locally free sheaf if and only if $E$ is nonspecial and globally generated.
\end{enumerate}
\end{lemma}
Lemma \ref{lem-CNYmain} transforms the weak Brill-Noether problem into showing that $\Phi_{X \to X}^{I_\Delta}(E)^{\vee}$ is a sheaf, and this is how Bridgeland stability becomes useful. Given the Mukai vector $(r, dH,a)$ of $E$,  there is a special chamber $\mathcal{C}$ in the $(s,t)$-plane right above the wall defined by 
$\mu_{\sigma}(I_x^\vee[1]) = \mu_{\sigma} (E)$, where $I_x$ is the ideal sheaf of a point $x \in X$.   Minamide, Yanagida and Yoshioka prove the following result.
\begin{theorem} \cite[Theorem 4.9]{MYY} 
If $\sigma \in \mathcal{C}$, then  the moduli space of Bridgeland stable objects  $M_{\sigma}(r, dH, a)$ is isomorphic to the Gieseker moduli space $M_{X, H}(a, dH, r)$ via $E\mapsto\Phi_{X \rightarrow X}^{I_{\Delta}}(E)^\vee$. 
\end{theorem}
Recall that for $t \gg 0$, the moduli space $M_{\sigma}(r,dH,a)$ is isomorphic to the Gieseker moduli space $M_{X, H}(r, dH, a)$.   As we decrease $t$ in order to reach the chamber $\mathcal{C}$, we may cross a number of Bridgeland walls.  If not all of the sheaves in $M_{X,H}(r,dH,a)$ get destabilized along the way,  then for $\sigma\in\mathcal{C}$, $M_{X, H}(r, dH, a)$ and $M_{\sigma}(r, dH, a)$ have a common open subset, so $M_{X, H}(r, dH, a)$ satisfies weak Brill-Noether by Lemma \ref{lem-CNYmain}. Furthermore, if the general sheaf in $M_{X, H}(a, dH, r)$ is locally free, then the general sheaf in $M_{X, H}(r, dH, a)$ is globally generated. 

However, along the way we may encounter totally semistable walls, where every sheaf gets destabilized. Bayer and Macr\`{i} \cite{BayerMacri, BayerMacri2} have classified the totally semistable walls. Using this classification, one can show that if ${\bf v} = (r, dH, a)$ is a Mukai vector with $r\geq 0$ ,$d >0$, and ${\bf v}^2\geq -2$ such that $M_{X, H}({\bf v})$ does not satisfy weak Brill-Noether, then there must exists a spherical Mukai vector ${\bf v}_1 = (r_1, d_1 H, a_1)$ satisfying
\begin{enumerate}
\item $0 < d_1 \leq d$
\item $0 > \langle {\bf v}, {\bf v}_1\rangle= 2ndd_1 -a r_1 - ra_1$
\item $0 < dr_1 -d_1r \leq da_1 - d_1a$.
\end{enumerate}
Determining whether there are totally semistable walls between the Gieseker chamber and $\mathcal{C}$ is a purely numerical problem, albeit a complicated one. 
 Moreover, we obtain a dichotomy. Either there are no totally semistable walls above the chamber $\mathcal{C}$ and weak Brill-Noether holds, or the largest totally semistable wall provides a resolution of the general sheaf in $M_{X,H}(r,dH,a)$ which can then be used to compute the cohomology. This strategy yields the following sharp qualitative theorem.

\begin{theorem}\cite{CoskunNuerYoshioka}
Let $X$ be a K3 surface such that $\Pic(X)\cong \mathbb{Z}H$ with $H^2=2n$. Let ${\bf v}=(r,dH,a)$ be a Mukai vector with ${\bf v}^2\geq -2$, $r\geq 2$ and $d>0$.
\begin{enumerate}
\item For each $r\geq 2$, there exists a finite set of tuples $(n,r,d,a)$ for which ${\bf v}$ fails to satisfy weak Brill-Noether. 
\item If $n\geq r$, then ${\bf v}$ satisfies weak Brill-Noether.
\item If $a \leq 1$, then ${\bf v}$ satisfies weak Brill-Noether.
    \item  If $d \geq r \left\lfloor \frac{r}{n}\right\rfloor +2$, then ${\bf v}$ satisfies weak Brill-Noether.
   \item Assume $a \geq 2$ and $n>1$.  If $n \geq 2r$ or $d\geq \left\lfloor \frac{2r}{n} \right\rfloor + 2$, then the general sheaf in $M_H({\bf v})$ is globally generated.
    \end{enumerate}
\end{theorem}

More importantly, the technique allows one to compute the cohomology of the general sheaf. For example, in \cite{CoskunNuerYoshioka} the authors classify all the Mukai vectors of rank at most 20 on K3 surfaces of Picard rank 1 for which weak Brill-Noether fails and compute the cohomology of the general sheaf in these cases. Recently, Liu has developed these ideas further and obtained an algorithm for computing the cohomology of spherical bundles on K3 surfaces \cite{Liuspherical}.

Similar strategies can be applied to other $K$-trivial surfaces such as abelian surfaces and Enriques surfaces. For example, in the case of abelian surfaces of Picard rank 1, one obtains the following.

\begin{theorem}\cite[Theorem 4.1]{CoskunNuer}
Let $X$ be an abelian surface with $\Pic(X)=\ZZ H$ and let ${\bf v} = (r, dH, a)$ be a Mukai vector such that $r \geq 0$, $d>0$ and ${\bf v}^2 \geq 6$. Then weak Brill-Noether holds for $M_{X,H}({\bf v})$.
\end{theorem}

\subsection{Asymptotic weak Brill-Noether} In general, the weak Brill-Noether problem is wide open. Given a rank $r \geq 2$ and a slope $\nu$, if $\Delta \gg 0$, then $\chi({\bf v}) < 0$. Moreover, by  O'Grady's Theorem \ref{thm-OGrady}, we can assume that the moduli space $M_{X, H}({\bf v})$ is irreducible and the general member is a slope-stable vector bundle. Potentially increasing $\Delta$, after a series of elementary modifications, we can assume that the general sheaf has no global sections. By applying the same argument to the Serre dual, we can conclude that the general sheaf also has no $h^2$.  We thus obtain an asymptotic weak Brill-Noether result.

\begin{proposition}\cite[Proposition 7]{CoskunHuizengaAbel}
Let ${\bf v}$ be a Chern character of rank at least $2$. If $\Delta \gg 0$, depending on $X, H, r$ and $\nu({\bf v})$, then $H^i(X, \cV) = 0$ unless $i=1$.
\end{proposition}

Similarly, by applying Serre Vanishing, one can show that the general sheaf in $M_{X,H}({\bf v} \otimes mA)$ has no higher cohomology if $A$ is ample and $m \gg 0$ (see \cite[Theorem 3.7]{CoskunHuizengaAbel}). It would be very useful to have good explicit bounds on $m$, especially for surfaces of general type.

\section{Applications and generalizations of weak Brill-Noether}\label{sec-applications}
In this section, we discuss some applications and generalizations of weak Brill-Noether Theorems.

\subsection{Classification of stable Chern characters}
The weak Brill-Noether problem is  closely tied to Problem \ref{prob-classifystabch}, the problem of classifying Chern characters of stable sheaves.

\subsubsection{The classification of stable vector bundles on $\PP^2$} Dr\'{e}zet and Le Potier classified Chern characters of stable bundles on $\PP^2$ \cite{DrezetLePotier, LePotier}. We now briefly recall this classification. 

\noindent{\bf Exceptional bundles.}  A coherent sheaf $\cV$ is  {\em exceptional} if $\Hom(\cV,\cV)=\CC$ and $\Ext^i(\cV,\cV)=0$ for $i > 0$.  Exceptional sheaves on $\PP^2$ have been classified by Dr\'{e}zet \cite{Drezet2}.  Let  $\cV$ be an exceptional sheaf on $\PP^2$ of rank $r$. First, since $\Ext^1(\cV, \cV) =0$, every small deformation of $\cV$ has to be trivial. Consequently, $g^*(\cV) \cong \cV$ for every $g \in \PP GL(3)$. We conclude that $\cV$ cannot have any singularities and must be a vector bundle.  By Riemann-Roch, $$\chi(\cV, \cV) = 1 = r^2(P(0) - 2\Delta) = r^2- c_1^2  + 2r \ch_2.$$ Hence, $r$ and $c_1$ are relatively prime and $\Delta(\cV) = \frac{1}{2}\left(1 - \frac{1}{r^2} \right).$ Observe that the discriminant of an exceptional bundle is less than $\frac{1}{2}$ and that the slope determines the Chern character of an exceptional bundle. 

If $E\subset\PP^2$ is an elliptic curve embedded as a cubic, then applying $\Hom (\cV, -)$ to the exact sequence 
$$0 \to \cV(-3) \to \cV \to \cV|_E \to 0$$ and using $$\Hom(\cV, \cV(-3))= \Ext^2(\cV(-3), \cV(-3))^* =0 \ \mbox{and} \ \Ext^1(\cV, \cV(-3))= \Ext^1(\cV(-3), \cV(-3))^*=0,$$ we conclude that $\Hom(\cV|_E, \cV|_E) = \CC$. Hence, the restriction of  $\cV$ to $E$ is simple. Since a bundle on an elliptic curve is a direct sum of its Harder-Narasimhan factors, we conclude that  $\cV|_E$ is semistable. As $E$ was arbitrary, it follows that $\cV$ is semistable and thus stable since $\gcd(r,c_1)=1$.  Consequently, an exceptional bundle on $\PP^2$  is stable. Conversely, a stable bundle $\cV$ on $\PP^2$ with $\Delta(\cV) < \frac{1}{2}$ is  exceptional by Riemann-Roch.  Moreover, two exceptional bundles with the same slope are  isomorphic since $\chi(\cV_1, \cV_2) =1$ and hence there  is a nontrivial homomorphism between $\cV_1$ and $\cV_2$ which must be an isomorphism by stability.

The line bundles $\OO_{\PP^2}(n)$ are exceptional. An {\em exceptional collection} on $\PP^2$ is a triple of exceptional bundles $(\cV_1, \cV_2, \cV_3)$ such that $\Ext^i(\cV_k, \cV_j)=0$ for $1 \leq j < k \leq 3$ and all $i$. On $\PP^2$, $(\OO_{\PP^2}, \OO_{\PP^2}(1), \OO_{\PP^2}(2))$ is the {\em standard exceptional collection}. Given a pair of adjacent bundles in an exceptional collection $\cV, \cW$, we can form the  left and right mutations of the pair
$$0 \to \cV \to \cW \otimes \Hom(\cV, \cW)^* \to R_{\cW}(\cV) \to 0$$
$$0 \to L_{\cV}(\cW) \to \cV \otimes \Hom(\cV, \cW) \to \cW  \to 0.$$
For example, the Euler sequence exhibits $T_{\PP^2}$ as $R_{\OO_{\PP^2}(1)}(\OO_{\PP^2})$.  Dr\'{e}zet proves that on $\PP^2$ the right $R_{\cW}(\cV)$  and left $L_{\cV}(\cW)$ mutations are again exceptional bundles. In fact, given an exceptional collection, $(\cV_1, \cV_2, \cV_3)$, the collections $(\cV_2, R_{\cV_2}(\cV_1), \cV_3)$ and $(\cV_1, L_{\cV_2}(\cV_3), \cV_2)$ are again exceptional collections. This justifies the terminology right/left mutation. Furthermore, every exceptional  bundle on $\PP^2$ is obtained from line bundles by a sequence of mutations \cite{Drezet2}. One can systematically generate the slopes of all exceptional bundles on $\PP^2$ by formalizing the process of mutations to obtain an explicit one-to-one correspondence $\varepsilon: \ZZ[\frac{1}{2}] \to \mathcal{E}$ between dyadic integers and the exceptional slopes, defined inductively by 
$\varepsilon(n) = n$ for an integer $n$ and 
$$\varepsilon \left( \frac{2p+1}{2^{q+1}} \right) = \varepsilon \left( \frac{p}{2^q} \right) .  \varepsilon \left( \frac{p+1}{2^q} \right),$$ where $$\alpha . \beta = \frac{\alpha+\beta}{2} + \frac{\Delta_{\beta} - \Delta_{\alpha}}{3+\alpha-\beta}.$$

\begin{figure}[t]

\begin{center}
\scalebox{.8}{\setlength{\unitlength}{1in}
\begin{picture}(5.34,2.15)
 \put(0,0){\includegraphics[scale=.5,bb=0 0 10.69in 4.31in]{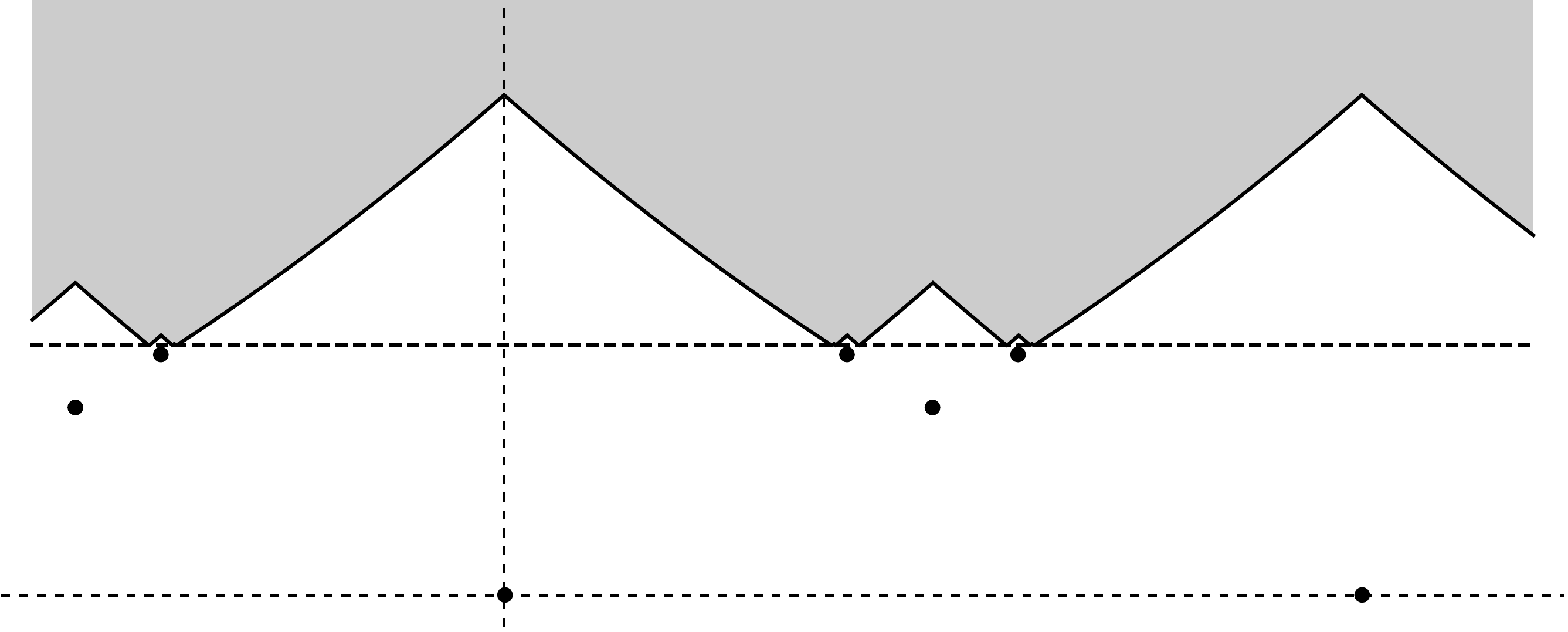}}
\put(5,.2){$\mu$}
\put(1.75,.2){$\OO_{\PP^2}$}
\put(4.2,.2){$\OO_{\PP^2}(1)$}
\put(2.9,.59){$T_{\PP^2}(-1)$}
\put(1.75,2){$\Delta$}
\put(1.75,1.8){$1$}
\put(1.75,1.05){$1/2$}
\put(1.75,.0){$0$}
\put(4.65,.0){$1$}
\put(3.5,1.45){$\Delta=\delta(\mu)$}
\end{picture}}
\end{center}
\caption{The Dr\'ezet--Le Potier curve.  Chern characters of stable bundles with positive dimensional moduli spaces lie in the shaded region above the curve.  The Chern characters of exceptional bundles are below the line $\Delta = \frac{1}{2}$. }\label{fig-DLP}
\end{figure}
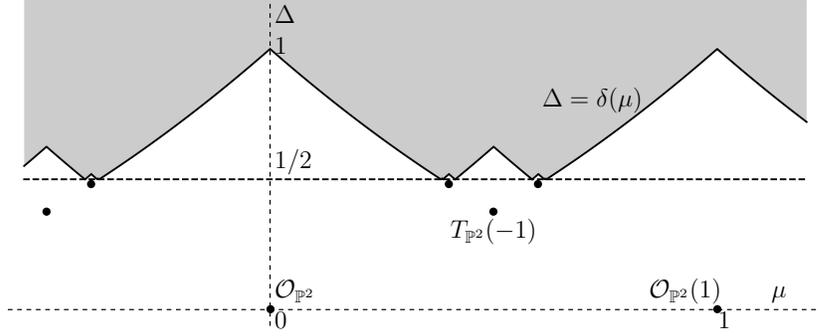

\noindent{\bf Stable bundles in general.} Given two stable bundles  on $\PP^2$ with $0 < \mu(\cV) - \mu(\cW)<3$, we have that $\Hom(\cV, \cW)=0$ and $\Ext^2(\cV, \cW) = \Hom(\cW, \cV(-3))=0$ by stability. Consequently, $$\chi(\cV, \cW) = r(\cV) r(\cW)\left(P \left(\mu(\cW) - \mu(\cV)\right) - \Delta(\cV) - \Delta(\cW)\right) \leq 0.$$ Setting $\cV$ to be an exceptional bundle $\mathcal{E}_{\alpha}$, gives an inequality for the discriminant $\Delta(\cW)$ of a stable bundle $\cW$ in terms of its slope $\mu(\cW)$ provided $0 < \mu(\mathcal{E}_{\alpha})- \mu(\cW) <3$. Similarly, setting $\cW$ to be an exceptional bundle $\mathcal{E}_{\alpha}$, gives an inequality for the discriminant $\Delta(\cV)$ of a stable bundle $\cV$ in terms of its slope $\mu(\cV)$ provided $0 < \mu(\cV) - \mu(\mathcal{E}_{\alpha}) <3$. Graphing the case of equality in the $(\mu, \Delta)$-plane for all  exceptional bundles yields a fractal curve called the  Dr\'ezet-Le Potier curve (see Figure \ref{fig-DLP}). The Chern character of any stable bundle which is not exceptional must lie above the Dr\'ezet-Le Potier curve depicted by the shaded region. Conversely, using dimension estimates, one can show that if the Chern character lies above the Dr\'ezet-Le Potier curve, then the general bundle defined by the Gaeta resolution is stable. One thus obtains the main classification theorem on $\PP^2$ due to Dr\'ezet and Le Potier.

\begin{theorem}[Dr\'ezet-Le Potier]
Let ${\bf v}$ be an integral Chern character of positive rank. There exists a Gieseker semistable sheaf with Chern character ${\bf v}$ if and only if $\Delta({\bf v}) \geq \delta(\mu({\bf v}))$ or ${\bf v}$ is a multiple of the Chern character of an exceptional bundle. When $\Delta({\bf v}) \geq \delta(\mu({\bf v}))$, the moduli space $M_{\PP^2, L}({\bf v})$ is an irreducible, normal projective variety of dimension $r^2(2\Delta -1) + 1$.  
\end{theorem}

\subsubsection{Classification of stable Chern characters on Hirzebruch surfaces} Using Theorems \ref{thm-HirzebruchWBN} and \ref{thm-Gaeta} as main tools, one can classify the Chern characters of stable sheaves on Hirzebruch surfaces  $\FF_e$ for generic polarizations (see \cite{CoskunHuizengaHExist}).

Let  $H_m = E + (m+e)F$ on $\FF_e$ with $m\in\QQ_{>0}$. Then $H_m$ is ample and every ample divisor is proportional to some $H_m$.  To classify Chern characters, one first shows that an $H_m$-semistable sheaf is $H_{\lceil m \rceil +1}$-prioritary and thus one needs to classify $H_k$-prioritary sheaves for positive integers $k$.  Then one computes the Harder-Narasimhan filtration of the general $H_k$-prioritary sheaf. The Gaeta-type resolution resolves the first issue. 

Let $\nu({\bf v}) = \epsilon E + \varphi F$. Set
$$\psi := \varphi + \frac{1}{2}e(\lceil \epsilon \rceil-\epsilon)-\frac{\Delta}{1-(\lceil \epsilon\rceil - \epsilon)},$$ and let $$L_0 := L_{\lceil \epsilon \rceil,\lceil \psi\rceil} = \lceil \epsilon\rceil E+\lceil \psi \rceil F.$$ Then the general sheaf in $\cP_{\FF_e, F} ({\bf v})$ admits a Gaeta-type resolution as in (\ref{eq-Gaeta}) with $L=L_0$. Expressing the condition that the sheaf be $H_k$ prioritary, one sees that we need $\chi({\bf v}(-L_0 -H_k))\leq 0$. Conversely, one can construct $H_k$-prioritary sheaves satisfying this inequality to conclude the following theorem.

\begin{theorem}\cite[Theorem 4.16]{CoskunHuizengaHExist}
Let ${\bf v}$ be an integral Chern character of positive rank on $\FF_e$ with $\Delta({\bf v}) \geq 0$ and let $k$ be a positive integer. Then the stack $\cP_{H_k}({\bf v})$ is nonempty if and only if
$$\chi ({\bf v}(-L_0-H_k)) \leq 0.$$
\end{theorem}
This theorem already provides stronger Bogomolov inequalities for the existence of semistable sheaves on $\FF_e$.

Next, we compute the Harder-Narasimhan filtration of the general $H_m$-prioritary sheaf. There exist $H_m$-semistable sheaves if and only if the filtration is trivial. Suppose the general $H_m$-Harder-Narasimhan filtration has length $\ell$ and the graded pieces have Chern characters ${\bf v}_i = (r_i, \nu_i, \Delta_i).$ Then $$\sum_{i=1}^{\ell} {\bf v}_i = {\bf v}.$$ Furthermore, the moduli spaces $M_{\FF_e, H_m}({\bf v}_i)$ are nonempty since the graded pieces are semistable sheaves. The fact that the sheaves are $H_{\lceil m \rceil +1}$-prioritary and the Schatz-stratum corresponding to this Harder-Narasimhan filtration has codimension 0 lead to additional inequalities that place strong restrictions on ${\bf v}_i$.

First, the prioritary condition implies that the restriction of the general sheaf to a general rational curve in the class $H_{\lceil m \rceil}$ or $H_{\lfloor m \rfloor}$ has balanced splitting. This translates to the inequality that $$|(\nu_i-\nu)\cdot H_m| < 1.$$
Next, a dimension computation shows that for the Schatz-stratum to have codimension 0,  the following orthogonality relations hold
 $$\chi({\bf v}_i, {\bf v}_j) = 0 \ \ \mbox{for} \ \ i < j.$$ Hence, the slopes $\nu_i$ are restricted to a bounded region, and since the ranks $r_i$ are bounded there are only finitely many possibilities for the $\nu_i$.  Furthermore,  the orthogonality relations imply the discriminant $\Delta_i$ has to be the minimal possible discriminant of an $H_m$-semistable sheaf with rank $r_i$ and total slope $\nu_i$. Hence, there are finitely many possible ${\bf v}_i$ that can be the Chern characters of the graded pieces of the generic $H_m$-Harder-Narasimhan filtration.
 
 Conversely, if one can find Chern characters ${\bf v}_i$, $1 \leq i \leq \ell$, that satisfy these constraints, then the general Harder-Narasimhan filtration has factors with these Chern characters. 
 \begin{theorem}\cite[Theorem 5.3]{CoskunHuizengaHExist}
Let ${\bf v}$ be a Chern character such that $\cP_{H_{\lceil m \rceil}}({\bf v})$ is nonempty. Let ${\bf v}_1, \dots, {\bf v}_{\ell} \in K(\FF_e)$ be positive rank Chern characters satisfying the following properties:
 \begin{enumerate}
 \item $\sum_{i=1}^{\ell} {\bf v}_i = {\bf v}$,
 \item The reduced Hilbert polynomials $q_i$ of ${\bf v}_i$ are strictly decreasing $q_1 > \cdots > q_{\ell}$,
 \item $\mu_{H_m}({\bf v}_1) - \mu_{H_m}({\bf v}_{\ell}) \leq 1$,
 \item $\chi({\bf v}_i, {\bf v}_j) =0$ for $i<j$,
 \item $M_{\FF_e, H_m}({\bf v}_i)$ is nonempty,
 \end{enumerate}
Then the Harder-Narasimhan filtration of the general sheaf in $\cP_{H_{\lceil m \rceil}}({\bf v})$ has length $\ell$ and the factors have Chern characters ${\bf v}_i$. 
 \end{theorem}
Thus determining the Harder-Narasimhan filtration of a general $H_{\lceil m \rceil}$-prioritary sheaf becomes a finite computational problem. In particular,  one obtains an algorithm for classifying Chern characters of $H_m$-semistable sheaves. Using the fact that $K(\F_e) \cong \ZZ^4$, one can show $\ell \leq 4$ to further simplify the problem.  

\begin{remark}
On $\PP^2$, either the subbundle or the quotient bundle in the general Harder-Narasimhan filtration of a general prioritary bundle is exceptional. Hence, one can give explicit inequalities for the discriminants of semistable bundles without computing the semistable bundles of lower rank. On $\FF_e$, it may happen that neither the quotient nor the subbundle in a generic Harder-Narasimhan filtration of length $2$ is an exceptional bundle (see \cite{CoskunHuizengaHExist} for explicit examples). In general, it is not enough to consider the constraints given by exceptional bundles.
\end{remark}

\subsection{The tensor product problem} A natural generalization of the Brill-Noether problem which plays a fundamental role in  the birational geometry of moduli spaces of sheaves is the tensor product problem (see \cite{ABCH, Huizenga}). 

\begin{problem}[The tensor product problem]
Let $\mathcal{V} \in M_{X, H}({\bf v})$ and $\mathcal{W} \in M_{X, H}({\bf w})$ be two stable sheaves on $X$. Compute the cohomology of $\cV \otimes \cW$.
\end{problem}
More generally, one can ask for the cohomology of  the tensor product of Bridgeland stable objects in the derived category of $X$.  The problem is already interesting when $\cV$ and $\cW$ are general elements in their moduli spaces. When $\cW = \OO_X$, the problem reduces to computing the cohomology of $\cV$. When $\cV$ is general in its moduli, we recover the weak Brill-Noether problem.

When $\cW=I_Z$ is an ideal sheaf of points and $\cV$ is a line bundle, then the problem reduces to the classical interpolation problem asking for when points impose independent conditions on sections of a line bundle. More generally, when $\cV$ is a higher rank sheaf, the problem is the higher rank interpolation problem asking for when points impose independent conditions on sections of a sheaf.  In general, these problems are wide open, though some important special cases have been solved. For example, the higher rank interpolation problem has an explicit combinatorial solution for zero-dimensional monomial schemes on $\PP^2$ \cite[Theorem 1.4]{CoskunHuizengaMonomial}.

On $\PP^2$ there is an almost complete answer to the tensor product problem for general sheaves. Let ${\bf v}$ be a stable Chern character on $\PP^2$ that lies on or above the Dr\'ezet-Le Potier curve and let $\cV$ be a general stable sheaf in $M_{\PP^2, L}({\bf v})$. Then the {\em associated exceptional bundle} $E_+$ is the exceptional bundle on $\PP^2$ with the smallest slope such that if $E'$ is any exceptional bundle with $\mu_L(E') > \mu_L(E_+)$, then $\chi(\cV \otimes E') >0$. It may happen that $\chi(\cV \otimes E_+) <0$, $\chi(\cV \otimes E_+) =0$ or $\chi(\cV \otimes E_+) >0$. By Theorem \cite[Theorem 4.1]{CoskunHuizengaWoolf}, the associated exceptional bundle exists. Briefly, by Riemann-Roch, the locus ${\bf v}^\perp := \{{\bf w} : \chi({\bf v} \otimes {\bf w}) = 0\}$ is an upward parabola in the slope-discriminant plane. This parabola intersects the line $\Delta=\frac{1}{2}$ in two points which lie under peaks of the Dr\'ezet-Le Potier curve. The Dr\'ezet-Le Potier curve has infinitely many peaks (see Figure \ref{fig-DLP}), each determined by an exceptional bundle.  The associated exceptional bundle $E_+$ is the exceptional bundle lying under the rightmost peak determined in this way.  (The points on the line $\Delta = \frac{1}{2}$ which do not lie under a peak form a generalized Cantor set $\mathfrak{C}$.  The fact that the parabola ${\bf v}^\perp$ meets $\Delta = \frac{1}{2}$ at points that are not in $\mathfrak{C}$ relies on number-theoretic properties of $\mathfrak{C}$.  For instance, the points of $\mathfrak{C}$ which are not endpoints of peaks have transcendental $\mu$-coordinate.)  

The answer to the problem is easiest to state in the case where $\chi({\bf v}\otimes E_+)\leq 0$.  In the other case where $\chi(\cV \otimes E_+) >0$, we need to fix some additional notation.  In this case we let $${\bf k} = {\bf v} - \chi(\cV\otimes E_+) \ch(E_+)$$ be the Chern character of the mapping cone of the canonical evaluation map $$E_+^* \otimes \Hom(E_+^*, \cV) \to \cV.$$

\begin{theorem}\cite[Theorem 1.2]{CoskunHuizengaKopper}
Let ${\bf v}$ and ${\bf w}$ be Chern characters of stable bundles on $\PP^2$. Suppose $M_{\PP^2, L}({\bf v})$ is positive dimensional and ${\bf w}$ is sufficiently divisible (depending on ${\bf v}$). Let $\cV \in M_{\PP^2, L}({\bf v})$ and $\cW \in M_{\PP^2, L}({\bf w})$ be general bundles.
\begin{enumerate}
\item If $\chi({\bf v} \otimes E_+) \leq 0$, then either $H^0(\cV \otimes \cW)=0$ or $H^1(\cV \otimes \cW)=0$.
\item If $\chi({\bf v} \otimes E_+) > 0$ and $\rk({\bf k}) \leq 0$,  then either $H^0(\cV \otimes \cW)=0$ or $H^1(\cV \otimes \cW)=0$.
\item Suppose $\chi({\bf v} \otimes E_+) > 0$ and $\rk({\bf k}) > 0$.
\begin{enumerate}
\item If $\chi({\bf k} \otimes {\bf w}) \geq 0$ or $\chi({\bf w} \otimes E_+^*) \leq 0$, then either $H^0(\cV \otimes \cW)=0$ or $H^1(\cV \otimes \cW)=0$.
\item Otherwise $\cV \otimes \cW$ is special and $$h^0(\cV \otimes \cW) = \chi({\bf v} \otimes E_+) \chi({\bf w} \otimes E_+^*) \quad \mbox{and} \quad h^1(\cV \otimes \cW) = - \chi({\bf k} \otimes {\bf w})$$
\end{enumerate}
\end{enumerate}
\end{theorem}

By applying Serre duality, the theorem computes the tensor product of two general stable bundles subject to the requirement that ${\bf w}$ is sufficiently divisible. One can give an explicit bound on how divisible ${\bf w}$  needs to be depending on ${\bf v}$. We do not know whether the divisibility is necessary on $\PP^2$ or an artifact of the proof. The following example shows that on $\PP^n$ for $n > 2$, some divisibility is necessary.

\begin{example}\cite[Example 1.3]{CoskunHuizengaSmith}
Let $\cV_m$ be the kernel of a general map 
$$0 \to \cV_m \to \OO_{\PP^3}(2)^{\oplus 4m} \to \OO_{\PP^3}(4)^{\oplus m} \to 0.$$ When $m=1$,  $h^0(\PP^3, \cV_1) = 6$ and $h^1(\PP^3, \cV_1) =1$. On the other hand, when $m \geq 2$, then $h^1(\PP^3, \cV_m) =0$.
\end{example}

\subsection{Construction of Brill-Noether divisors and birational geometry}

Let ${\bf v}$ and ${\bf w}$ be two stable Chern characters on $X$ such that $\chi({\bf v} \otimes {\bf w}) = 0$. Then given a general sheaf $\cW \in M_{X, H} ({\bf w})$, we can define the virtual Brill-Noether divisor 
$$D_{\cW} := \{ \cV \in M_{X,H}({\bf v}) | h^1 (\cV \otimes \cW) \neq 0 \}.$$
In general, $D_{\cW}$ may not be a divisor. If $h^i(\cW \otimes \cV') =0$ for all $i$ and a general $\cV'$, then $D_{\cW}$ is an effective divisor. Hence, the tensor product problem in the special case of $\chi({\bf v} \otimes {\bf w}) = 0$ is the key ingredient for constructing effective Brill-Noether divisors.

There has been significant progress in computing the ample and effective cones of $M_{X, H}({\bf v})$ and running the minimal model program for these moduli spaces in recent years, especially when $X$ is $\PP^2$, a Hirzebruch, K3, abelian or Enriques surface. Surveying the developments in this direction would take us too far afield, so we refer the reader to  \cite{ABCH, BayerMacri, BayerMacri2, BertramCoskun, CoskunHuizengaAmpleP2, CoskunHuizengaGokova, CoskunHuizengaAmple, CoskunHuizengaWoolf, Huizenga, LiZhao, MYY, Nuer16a, NuerYoshioka} for some of these developments.

\subsection{Constructions of Ulrich bundles on surfaces} A solution of the weak Brill-Noether problem allows one to construct Ulrich bundles on surfaces. 

\begin{definition}
Let $X \subset \PP^n$ be a smooth, projective variety of dimension $d$. An {\em Ulrich bundle} $\mathcal{V}$ on $X$ is a bundle that satisfies $H^i(X, V(-j))=0$ for $1 \leq j \leq d$ and all $i$.
\end{definition}
Ulrich bundles play a central role in the study of Chow forms of a variety \cite{ESW},  the minimal resolution conjecture (see \cite{AGO}) and Boij-S\"{o}derberg Theory  (see \cite{EisenbudSchreyer}). For example, Eisenbud and Schreyer show that the cone of  cohomology tables of $X$ is the same as that of $\PP^d$ if and only if $X$ admits an Ulrich bundle \cite{EisenbudSchreyer}. Eisenbud and Schreyer raise the question whether every  projective variety admits an Ulrich bundle. Existence is known in some cases including smooth curves \cite{ESW}, complete intersections \cite{BaHU},  del Pezzo surfaces \cite{CKM}, certain rational surfaces \cite{ESW, Kim}, K3 surfaces \cite{AGO, Faenzi2}, abelian surfaces \cite{BeauvilleU} and certain Enriques surfaces \cite{BorisovNuer} among many others. 

Let $\mathcal{V}$ be an Ulrich bundle for a polarized surface $(X, H)$ with rank $r$, total slope $\nu$ and discriminant $\Delta$.  Since $\mathcal{V}(-H)$ and $\mathcal{V}(-2H)$ have no cohomology, the Euler characteristics of $\mathcal{V}(-H)$ and $\mathcal{V}(-2H)$ must vanish. By Riemann-Roch, we conclude that
$$\frac{1}{2} \nu^2 - \frac{1}{2} \nu \cdot K_X + \chi(\OO_X) - \Delta =0$$
$$\frac{1}{2} (\nu-H)^2 - \frac{1}{2} (\nu-H) \cdot K_X + \chi(\OO_X) - \Delta =0$$
Consequently, an Ulrich bundle on a surface satisfies 
\begin{equation}\label{eq-Ulrich}
2 \nu \cdot H = H^2 + H \cdot K_X \quad \mbox{and} \quad 2 \Delta = \nu^2 - \nu \cdot K_X + 2 \chi(\OO_X).
\end{equation}
Conversely, if we take $2\nu = H + K_X$,  $M_{X,H}(r, \nu, \Delta)$ contains locally free sheaves, and it satisfies weak Brill-Noether, then the general locally free sheaf in  $M_{X,H}(r, \nu, \Delta)$  is an Ulrich bundle. 

Combining O'Grady's Theorem with Serre vanishing yields an asymptotic existence result on any smooth projective surface.

\begin{theorem}\cite[Theorem 4.3]{CoskunHuizengaICMSat}\label{thm-Ulrichasymptotic}
Let $(X, H)$ be a smooth, polarized surface. There exists a positive integer $m_0$ such that for all $m \geq m_0$, the polarized surface $(X, mH)$ admits an Ulrich bundle of every positive even rank. Moreover, if $K_X$ (respectively, $K_X + H)$ is divisible by $2$ in $\Pic(X)$ and $2m \geq m_0$ (respectively, $2m+1 \geq m_0$), then $(X, 2mH)$ (respectively, $(X, (2m+1)H)$) admits an Ulrich bundle of every rank $r \geq 2$.
\end{theorem}

A polarized variety $(X,H)$ is of {\em Ulrich wild representation type} if one can find arbitrarily large dimensional families of Ulrich bundles on $X$ with respect to $H$. An easy corollary Theorem \ref{thm-Ulrichasymptotic} is the following:

\begin{corollary}\cite[Corollary 4.5]{CoskunHuizengaICMSat}
Let $(X,H)$ be a smooth, polarized surface. Then there exists an integer $m_0$ such that for all $m \geq m_0$, $(X, mH)$ is of Ulrich wild representation type.
\end{corollary} 

On specific surfaces, weak Brill-Noether gives detailed information on Chern classes of Ulrich bundles. For example, on Hirzebruch surfaces one obtains the following (see also  V. Antonelli \cite{Antonelli}).  

\begin{theorem}\cite[Theorem 4.6]{CoskunHuizengaICMSat}
Let $H = aE+bF$ be an ample divisor on $\FF_e$. Let ${\bf v} = (r, \nu, \Delta) = (r, \alpha E + \beta F, \Delta)$ be an integral Chern character with $r\geq 2$.  There exists a locally free $F$-prioritary sheaf $\cV$ with Chern character ${\bf v}$ satisfying $$H^i(\FF_e, \cV(-H)) = H^i(\FF_e, \cV(-2H)) =0 \ \forall i$$ if and only if $$ a-1 + \frac{ea(a-1)}{2b} \leq \alpha \leq 2a-1-\frac{ea(a-1)}{2b},$$  $$\beta = \left( e- \frac{b}{a}\right) (\alpha +1) + 3b-1 - \frac{e}{2}(3a+1)$$ and $$\Delta(\cV) =  \left(\frac{e}{2} - \frac{b}{a} \right) (\alpha^2 + (2-3a)\alpha +2a^2 -3a+1).$$
In particular, there exists an Ulrich bundle of rank $2$ for $H= aE + bF$ with invariants $${\bf v} = (r, \nu, \Delta) = \left(2, \left(\frac{3}{2}a-1\right)E + \left(\frac{3}{2}b - \frac{e}{2}-1\right)F, \frac{a}{8}(2b-ae)\right).$$
\end{theorem}

Similarly detailed theorems can be proven on certain blowups of $\PP^2$ or $\FF_e$ using the corresponding weak Brill-Noether statements. We refer the reader to \cite[Theorem 4.10 and 4.11]{CoskunHuizengaICMSat} for precise statements. The weak Brill-Noether Theorems for K3 or abelian surfaces similarly yield classification results for Ulrich bundles. We state the result for K3 surfaces and refer the reader to \cite[Corollary 6.1]{CoskunNuer} for abelian surfaces.

\begin{proposition}\cite[Proposition 4.4]{CoskunNuerYoshioka}
Let $X$ be a K3 surface with $\Pic(X)=\ZZ H$.  There exists an Ulrich bundle of rank $r$ with respect to $mH$ if and only if $2\mid rm$.  Moreover, when an Ulrich bundle of rank $r$ exists, it has Mukai vector ${\bf v}=\left(r,\left(\frac{3rm}{2}\right)H,r(2m^2n-1)\right)$.  In particular, there exists an Ulrich bundle of any rank $r\geq 2$ with respect to $2H$.
\end{proposition}

One can also prove results on variants of Ulrich bundles. For example, Eisenbud and Schreyer \cite{EisenbudSchreyer} conjectured that there exist bundles $\cV$ on $\PP^1 \times \PP^1$ with natural cohomology, meaning that $\cV(a F_1 + bF_2)$ has at most one nonzero cohomology group for any $a, b \in \ZZ$. The weak Brill-Noether theorem on $\PP^1 \times \PP^1$ immediately yields a positive solution (see also \cite{Solis}).

\begin{corollary}\cite[Corollary 4.13]{CoskunHuizengaICMSat}
Let $F_1$ and $F_2$ denote the two fiber classes on $\PP^1 \times \PP^1$. 
Let ${\bf v}$ be an integral Chern character such that $\rk({\bf v}) \geq 2$ and $\Delta \geq 0$. Then $\cP_{\PP^1 \times \PP^1, F_1}({\bf v})$ is nonempty, and the general $\cV \in \cP_{\PP^1 \times \PP^1, F_1}({\bf v})$ is locally free and has at most one nonzero cohomology group. In particular, the very general member of $\cP_{\PP^1 \times \PP^1, F_1}({\bf v})$ is a bundle with natural cohomology.
\end{corollary}

\section{Cohomology jumping loci}\label{sec-jumping}
When the weak Brill-Noether Theorem holds for  $M_{X, H}({\bf v})$, it is natural to try to describe the loci of bundles with unexpected cohomology. In this section, we will review recent developments on the Brill-Noether loci on $\PP^2$ due to Gould, Lee and Liu \cite{GouldLiu} and describe some generalizations to other surfaces.

Let $M_{X,H}({\bf v})$ be an irreducible moduli space such that every sheaf $\cV$ has $H^2(X, \cV)=0$. For example, this condition automatically holds on moduli spaces of sheaves on $\PP^2$, a K3, an abelian or a Hirzebruch surface if $\nu({\bf v})$ is effective. Assume that the locus of stable sheaves $M_{X,H}({\bf v})^s$ is nonempty. Define the $k$-th Brill-Noether locus $B_{X,H}^k({\bf v})$ as the closure of the locus of stable sheaves with at least $k$-independent sections:
$$B_{X,H}^k({\bf v}) := \overline{\{ \cV \in M_{X,H}({\bf v})^s | h^0(X, \cV) \geq k\} }\subset M_{X, H}({\bf v}).$$
The locus $B_{X,H}^k({\bf v}) \cap M_{X,H}({\bf v})^s$ has a natural determinantal structure (see \cite[Proposition 2.7]{CoskunHuizengaWoolf} or \cite[\S 2]{CostaMiroRoig}), so every irreducible component of $B_{X,H}^k({\bf v})$ has codimension at most  the expected value $k (k- \chi({\bf v}))$.  By definition, we have $B_{X,H}^{k+1}({\bf v}) \subset B_{X, H}^k({\bf v})$ for $k \geq 0$.

\begin{problem}
One can immediately raise the following natural problems.
\begin{enumerate}
\item Determine when $B_{X, H}^k({\bf v})$ is nonempty, equivalently determine $\max_{\cV \in M_{X,H}({\bf v})^s} (h^0(X,\cV))$.
\item Classify the irreducible components of  $B_{X, H}^k({\bf v})$  and determine their dimensions. In particular, determine when $B_{X,H}^k({\bf v})$ is irreducible.
\item Describe the singularities of $B_{X, H}^k({\bf v})$.  
\item When do Torelli type theorems hold for $B_{X,H}^k({\bf v})$, i.e., when can  $X$  be recovered from $B_{X,H}^k({\bf v})$?
\end{enumerate}
\end{problem}
In general, very little is known about these problems. Even in rank 1, the loci $B_{X,H}^k({\bf v})$ can have many components of different dimensions and are poorly understood.

\begin{example}\label{ex-redrank1}  \cite[Proposition 5.3]{GouldLiu} Let ${\bf v} = \ch(I_Z(3))$ on $\PP^2$, where $Z$ is a zero-dimensional scheme of length $|Z| = n \geq 10$. The Brill-Noether locus $B_{\PP^2, L}^2({\bf v})$ parameterizing $I_Z(3)$ with $h^0(\PP^2, I_Z(3)) \geq 2$  has at least two irreducible components. 
\begin{enumerate}
\item Let $W_1 \subset \PP^{2[n]}$ be the closure of the locus  with  $n-4$ points lying on a  line and 4 general points. Then $\dim(W_1) = n+6$. 
\item Let $W_2 \subset \PP^{2[n]}$ be the closure of the  locus with $n-1$ points lying on a  smooth conic and 1 general point. Then $\dim(W_2) =n+6$. 
\end{enumerate}
We claim that if  $I_Z(3) \in B_{\PP^2,L}^2({\bf v})$ where $Z$ has distinct points, then $Z \in W_1$ or $Z \in W_2$. Let $C_1$ and $C_2$ be two distinct curves of degree 3 vanishing on $Z$. By B\'ezout's Theorem,  $C_1$ and $C_2$ must have a common curve $D$. Then $D$ is either a conic  (possibly reducible or nonreduced) or a line.  If $D$ is a conic, then the residual curves must be lines and they intersect in at most one point. Hence,  $Z \in W_2$. If $D$ is a line, then the residual curves must be conics. Since they do not have a common curve, their intersection must be a scheme of length 4. Hence, $Z \in W_1$. Since $W_1$ and $W_2$ have the same dimension, they form distinct irreducible components of $B_{\PP^2, L}^2({\bf v})$. 

A slightly more subtle argument allows Gould, Lee and Liu to deduce that when ${\bf v} = \ch(I_Z(k))$, where $|Z| > k^2$, the Brill-Noether locus $B_{\PP^2, L}^2({\bf v})$ has at least $k$ components, typically not all of the same dimension.
\end{example}

\begin{example}\label{ex-rank1general}
The previous example is not special to $\PP^2$.  For simplicity, let $X$ be a surface of Picard rank 1 generated by an effective ample divisor $H$ and assume  $H^1(X, \OO_X) =0$.  Let $m, n$ be two integers such that $m> 2$ is a sufficiently large odd integer and $n > m^2 H^2$. Let  ${\bf v} = \ch(I_Z(mH))$ with $|Z|=n$.  For $0 < s < m$, let $W_s \subset X^{[n]}$ be the closure of the locus of points with $h^0(X, \OO_X((m-s)H))-2$ general points on $X$ and $n+ 2 - h^0(X, \OO_X((m-s)H))$ general points on a  curve with class $sH$. Then
$$\dim(W_s) = h^0(X, \OO_X(sH)) + h^0(X, \OO_X((m-s)H)) + n-2.$$ Note that $\dim(W_s) = \dim(W_{m-s})$ and if $s \gg 0$, then we can assume that the curve of degree $s$ is irreducible. Hence, there are at least two components $W_{s_0}$ and $W_{m-s_0}$ with maximal dimension.  The inclusion $I_Z(mH) \subset \OO_X(mH)$ induces two independent curves $C_1$ and $C_2$ in $|mH|$ that contain $Z$. Let $I_Z(mH)$ be a general point of an irreducible component of $B_{X,H}^2({\bf v})$. Since $|Z| > m^2 H^2$, by B\'ezout's Theorem, $C_1$ and $C_2$ must have a common component of class $kH$. By the semicontinuity of the degree of the base curve, this component cannot contain $W_s$ for $s > k$. Hence, if $B_{X,H}^2({\bf v})$ is irreducible, we must have $k=1$. There can be at most $h^0(X, \OO_X((m-1)H))-2$ general points which do not lie on the common component. Hence, if this locus contains $W_1$, it must coincide with $W_1$. Then at least one of $W_{s_0}$ or $W_{m-s_0}$ is a component distinct from $W_1$ and we conclude that $B_{X,H}^2({\bf v})$ is reducible. 

In fact, by the semicontinuity of the degree of the base locus, observe that $W_s \not\subset W_t$ if $s > t$.  Assuming that there are irreducible sections of $H^0(X, \OO_X(kH))$ for every $k \geq 1$, then the loci $W_s$ have no containments among each other and one can show that $B^2_{X,H}({\bf v})$ has at least $k$ components.
\end{example}

Once the Brill-Noether loci in rank 1 are reducible, one can make examples of higher rank reducible Brill-Noether loci.

\begin{example}\cite[Example 5.5]{GouldLiu}
Consider extensions of the form 
$$0 \to \OO_{\PP^2} \to \cV \to I_Z(3) \to 0,$$ where $|Z| \geq 10$. Then $\chi(\cV) = 11- |Z|$. The Cayley-Bacharach condition is trivially satisfied, so there exists vector bundles of this form. Furthermore,  if the points of $Z$ are not collinear then there are no nonzero maps $\OO_{\PP^2} (2) \to  I_Z(3)$, so $\cV$ is stable. We have $$\ext^1(I_Z(3), \OO_{\PP^2}) = |Z|-1.$$ Set ${\bf v} = \ch(\cV)$ and ${\bf w} = \ch(I_Z(3))$. Now assume that $Z \in W_i$ for $i=1$ or 2 as in Example \ref{ex-redrank1}. Then $\cV \in B^3_{\PP^2, L} ({\bf v})$. We thus obtain two loci in $B^3_{\PP^2, L}({\bf v})$ of dimension $2|Z| +2$. The expected codimension of $B^3_{\PP^2, L} ({\bf v})$ is $3(|Z|-8)$. Since the moduli space has dimension $4|Z|-12$, the expected dimension of $B^3_{\PP^2, L} ({\bf v})$ is $|Z|+12$. Observe that we obtain unexpectedly large dimensional loci in $B^3_{\PP^2, L}({\bf v})$ if $|Z| > 10$. One can also easily check that the two components of $B^2_{\PP^2, L} ({\bf w})$ yield two distinct components of $B^3_{\PP^2, L}({\bf v})$. 
\end{example}

Similar constructions using Example \ref{ex-rank1general} produce reducible Brill-Noether loci on more general surfaces. Subbundles with a large space of sections is another source of components of Brill-Noether loci.

\begin{example}\cite[Example 5.11]{GouldLiu}
Let ${\bf v}_k$ be a a Chern character on $\PP^2$ with $\mu_L({\bf v}_k) = \frac{2}{3}$ and $\Delta({\bf v}_k) = \frac{5}{9} + \frac{k}{3}$. Then one can consider the following two types of bundles in $B^3_{\PP^2, L}({\bf v}_k)$.
\begin{enumerate}
\item  Let  $C$ be a smooth conic, $D$ be a divisor of degree $1-k$ and $\cV$ be a general extension of the form
$$ 0 \to \OO_{\PP^2}^{\oplus 3} \to \cV \to \OO_C(D) \to 0.$$
\item Let $T_{\PP^2}$ be the tangent bundle of $\PP^2$,  $Z$ be a zero-dimensional scheme of length $k+1$ and $\cV$ be defined by a nonsplit extension of the form
$$0 \to T_{\PP^2}(-1) \to \cV \to I_Z(1) \to 0.$$ 
\end{enumerate}
One can check that in both cases $\cV$ is stable and give distinct irreducible components of $B^3_{\PP^2, L} ({\bf v}_k)$ when $k \gg 0$.
\end{example}

Using a systematic analysis of the ideas in these examples, Gould, Lee and Liu prove the following structure theorem on Brill-Noether loci in moduli spaces of sheaves on  $\PP^2$. 

\begin{theorem}\cite[Theorem 1.1]{GouldLiu}
Let $M_{\PP^2, L} ({\bf v})$ be a nonempty moduli space of sheaves on $\PP^2$ with $\mu_L({\bf v}) > 0$.
\begin{enumerate}
\item For any $\cV \in M_{\PP^2, L} ({\bf v})$, $$h^0(\PP^2, \cV) \leq \max\left( r(\cV), P(c_1(\cV))= \frac{1}{2} c_1(\cV)^2 + \frac{3}{2} c_1(\cV) + 1\right).$$ In particular, $B^k({\bf v}) = \emptyset$ if $k > \max(r({\bf v}), P(c_1({\bf v})))$.
\item The Brill-Noether locus $B^r({\bf v})$ is nonempty.
\item If $\ch_1({\bf v}) = L$, then all the nonempty Brill-Noether loci are irreducible of the expected dimension.
\item If $\mu_L({\bf v}) > \frac{1}{2}$ and is not an integer and $\Delta({\bf v}) \gg 0$, then $B^r({\bf v})$ is reducible and contains components of both the expected and larger than expected dimension.
\end{enumerate}
\end{theorem}

The irreducibility in part (3) of the theorem also holds for more general surfaces.  Let $(X, H)$ be a polarized surface.  Let ${\bf v} = (r, C, d)$ be a Chern character on $X$ with $C \cdot H >0$. Assume that there does not exist a class $D$ on $X$ with $0 <D\cdot H < C \cdot H$. On $\PP^2$, the class of a line $L$ satisfies this property. More generally, if $X$ has Picard rank $1$, then the ample generator of the Picard group has this property. In this case, the Brill-Noether loci are better behaved.  The linear system $|C|$ cannot have any reducible or nonreduced elements (it can however be empty if $C$ is not effective). Given a degree $\kappa$, let $$\gamma_{\kappa} := \max\{h^0(B, \cL) | B \in |C|, \deg(\cL) = \kappa\}.$$ We have the easy bound $\gamma_{\kappa} \leq \kappa+1$, hence $\gamma_\kappa =0$ when $\kappa<0$. Combining the Riemann-Roch Theorem with Clifford's Theorem, one can obtain  better bounds on $\gamma_{\kappa}$.

\begin{proposition}\label{prop-mainpropBN}
Let $(X, H)$ be a polarized surface. Let ${\bf v}= (r, C, d)$ be a Chern character on $X$ with $C \cdot H >0$. Assume that there does not exists a divisor class $D$ with $0 < D \cdot H < C \cdot H$. Set $\kappa= d - \frac{C^2}{2}$. Then 
$$h^0(X, \cV) \leq r + \gamma_\kappa.$$
In particular, if $\kappa < 0$, then $h^0(\mathcal{V}) \leq r$ for every $\cV \in M_{X, H}({\bf v})$ and $B^k_{X, H}({\bf v}) = \emptyset$ for $k>r$.
Conversely, if $C$ is effective and $\Delta({\bf v}) \gg 0$ (equivalently, $d \ll 0$), then $B^{r}_{X, H}({\bf v})$ is nonempty.
\end{proposition}

\begin{proof}
Let $\cV \in M_{X, H}({\bf v})$ be a stable sheaf such that  $h^0(X, \cV) \geq r$. Choosing $r$ linearly independent sections induces a map $\phi: \OO_X^{\oplus r} \to \cV$. We claim that $\phi$ is an injective sheaf map. Suppose the image of $\phi$ is a sheaf $\cW$ of rank $s < r$. Since $h^0(X, \cW) \geq r$, we conclude that $\mu_H(\cW) > 0$. Since $C$ is the minimal class with $0 < C \cdot H$, we have that $c_1(\cW) \cdot H \geq C \cdot H$. Since $s < r$, this contradicts the stability of $\cV$. We conclude that  the image of $\phi$ must have rank $r$, hence its kernel is zero. So we have a short exact sequence
\begin{equation}\label{eq-extBN}
0 \to \OO_X^{\oplus r} \to \cV \to \OO_{C}(D) \to 0,
\end{equation}
where $\OO_C(D)$ is a pure rank one sheaf on a curve of class $C$ since the quotient of a torsion-free sheaf by a locally free sheaf cannot have torsion supported in dimension zero. \cite[Example 1.1.16]{HuybrechtsLehn}. 
By the long exact sequence on cohomology we conclude that $$h^0(X, \cV) \leq r + h^0(X, \OO_C(D)) \leq r + \gamma_{\kappa}.$$
In particular, if $\kappa <  0$, then $h^0(X, \cV) \leq r$. 

If $C$ is effective and $\Delta \gg 0$ (or equivalently, $d \ll 0$), then there exists nonsplit extensions of the form (\ref{eq-extBN}). It suffices to take $\Delta$ large enough so that $\chi(\OO_C(D), \OO_X) < 0$. Note that $\cV$ is torsion-free  (see \cite[Lemma 4.5]{GouldLiu}), stable and $h^0(X, \cV) \geq r$. Hence, $B^{r}_{X, H}({\bf v})$ is nonempty.
\end{proof}

\begin{remark}
When $c_1(\cV)$ is not minimal, one can use the Grauert-M\"{u}lich Theorem and restriction to curves to give a bound on $\max_{\cV \in M_{X, H}({\bf v})^s}(h^0(X, \cV))$ (see \cite[\S 3.3]{HuybrechtsLehn}).
\end{remark}

\begin{theorem}\label{thm-BNmain}
Let $(X, H)$ be a polarized surface such that $K_X \cdot H \leq 0$. Let ${\bf v}= (r, C, d)$ be a Chern character on $X$ with $C \cdot H > 0$ and $\Delta({\bf v}) \gg 0$. Assume that $h^0(X, \cL) > 0$ and $H^i(X, \cL) =0$ for $i \geq 1$ for every $\cL$ with $c_1(\cL) = C$. Further assume that there does not exist a divisor class $D$ with $0 < D \cdot H < C \cdot H$.  Then $B^k_{X, H} ({\bf v})$ is empty if $k > r$ and nonempty, irreducible and of the expected dimension if $k \leq r$.
\end{theorem}

\begin{proof}
We prove the theorem by induction on the rank. We already know that $B^k_{X, H} ({\bf v}) = \emptyset$ if $k > r$ by Proposition \ref{prop-mainpropBN}. It is also easy to see this when $r=1$. An element of $M_{X, H}({\bf v})$ is of the form $I_Z(\cL)$, where $\cL$ is a line bundle with class $C$. Since $\Delta({\bf v}) \gg 0$, we have $|Z| \gg 0$. We only need that $|Z| > C^2$. By B\'ezout's Theorem, if there exists two independent curves $C_1, C_2$ in $|\cL|$ containing $Z$, then they must have a common component. This is a contradiction since there are no reducible or nonreduced curves in the linear system $|\cL|$ by the minimality of $C$. 

On the other hand, if $k=1$, then $Z$ is a zero-dimensional scheme contained in an element of some linear system $|\cL|$. Since every curve in $|\cL|$ is reduced and irreducible and $X$ is smooth, the Hilbert scheme on the curve is irreducible of dimension $|Z|$. We conclude that $B^1_{X, H} ({\bf v})$ is irreducible of dimension $$|Z| + h^0(X, C) + \dim(\Pic^0(X)) -1 = \dim(M_{X, H} ({\bf v}) - 1+ \chi({\bf v})$$ as expected. 

Suppose that the theorem holds for rank less than $r$. By Proposition \ref{prop-mainpropBN}, we may assume that $k\leq r$. For $\cV \in B^k_{X, H} ({\bf v})$, consider the evaluation map $\phi: \OO_X^{\oplus k} \to \cV$. Then the saturation $\cW$  of the image of $\phi$ has rank $s \leq k$. If $s < k$, then $\cW$ is a subsheaf of $\cV$ with $\mu_H (\cW) > 0$. Since $C \cdot H$ is the minimal possible positive value and $r(\cW) < r$, we get a contradiction to the stability of $\cV$. We conclude that the evaluation map is injective. Hence, we obtain an exact sequence
$$0 \to \OO_X^{\oplus k} \to \cV \to \cV' \to 0.$$ 
As above, $\cV'$ cannot have torsion supported in dimension zero as the quotient of a torsion-free sheaf by a locally free one. In fact, when $k<r$, $\cV'$ is a $\mu_H$-stable sheaf.   Indeed, if $\cV'$ had either one dimensional torsion or a destabilizing subsheaf, then since the $H$-degree of $\cV'$ is the minimal, taking the quotient by the torsion or destabilizing subsheaf would yield a torsion free sheaf $Q$ with $\mu_H (Q) \leq 0$.  Since $Q$ would also be a quotient of $\cV$, this would contradict the stability of $\cV$.  So we conclude that when $k<r$ any $\cV\in B^k_{X, H} ({\bf v})$ is a nonsplit extension of a $\mu_H$-stable $\cV'$ of Chern character ${\bf v}-k{\bf v}(\OO_X)$ by $\OO_X^{\oplus k}$.

Conversely, consider nonsplit extensions of the form 
$$0 \to \OO_X^{\oplus k} \to \cV \to \cV' \to 0.$$
We have $\hom(\cV', \OO_X) =0$ by stability. 
We split the proof into two cases, $K_X \cdot H < 0$ and $K_X\cdot H=0$.

Suppose first that $K_X\cdot H<0$.  If $k < r-1$, then $\ext^2(\cV', \OO_X) = \hom(\OO_X, \cV'(K_X))=0$. The latter follows from the fact that $\mu_H(\cV'(K_X))<0$ since $r(\cV') \geq 2$ and $c_1(\cV')$ has minimal positive possible $H$-degree. If $k = r-1$ and $K_X \cdot H < 0$, then $\ext^2(\cV', \OO_X) = \hom(\OO_X, \cV'(K_X))=0$ unless $\cV'(K_X)=\OO_X$. The latter cannot happen since $\Delta({\bf v}) \gg 0$. We conclude that $\ext^1(\cV', \OO_X^{\oplus k})$ is constant for every $\cV' \in M_{X, H}({\bf v}')$ of dimension $- k \chi(\cV',  \OO_X)$. Since $\Delta({\bf v}) \gg 0$, we must also have $\Delta({\bf v}') \gg 0$. By O'Grady's Theorem, the moduli space  $M_{X, H}({\bf v}')$ is irreducible of the expected dimension. We conclude that $B^k_{X, H} ({\bf v})$ is the image of a projective bundle over an irreducible space. Hence, $B^k_{X, H}({\bf v})$ is irreducible and a dimension count shows that 
$$\dim(B^k_{X, H } ({\bf v}))= \dim (M_{X, H}({\bf v}')) - k \chi(\cV', \OO_X) - k^2.$$
We have  ${\bf v}' = (r-k, c, d)$. Since $\Delta({\bf v}), \Delta({\bf v}') \gg 0$, we have that 
$$\dim(M_{X, H}({\bf v})) = 2 r^2 \Delta({\bf v}) - (r^2 -1) \chi(\OO_X) \ \mbox{and} \ \dim(M_{X, H}({\bf v}')) = 2 r^2 \Delta({\bf v}') - ((r-k)^2 -1) \chi(\OO_X)$$ by O'Grady's Theorem. Hence, by Riemann-Roch
$$\dim(B^k_{X, H} ({\bf v})) = C^2 + (k-2d)r + (kr-r^2 +1) \chi(\OO_X) - \frac{k}{2} C \cdot K_X  - k^2.$$ 
On the other hand, the expected dimension of $B^k_{X, H} ({\bf v})$ is  given by 
$$\dim(M_{X, H} ({\bf v})) - k (k - \chi({\bf v}))= 2r^2 \Delta({\bf v}) - (r^2 -1) \chi(\OO_X)  - k (k - \chi({\bf v}))$$ which by Riemann-Roch is equal to
$$ C^2 + (k-2d)r + (kr-r^2 +1) \chi(\OO_X) - \frac{k}{2} C \cdot K_X  - k^2.$$ Hence, $B^k_{X, H}({\bf v})$ has the expected dimension. 

When $k=r$, then we get an extension of the form (\ref{eq-extBN}) and the irreducibility follows by a similar argument using the irreducibility of the compactified Jacobian. 

Finally, if $K_X \cdot H =0$, then $\hom(\OO_X, \cV'(K_X)) = s$ along the Brill-Noether locus $B^s_{X, H} ({\bf v}')$ which by induction on the rank has codimension $s (s- \chi({\bf v}' \otimes K_X))$. Since the latter is very large when $\Delta \gg 0$, these loci do not contribute new components. We conclude that the same estimate holds in this case as well. 
We conclude that the Brill-Noether locus $B^k_{X, H } ({\bf v})$ is irreducible, nonempty and of the expected dimension if $k \leq r$ and $\Delta({\bf v}) \gg 0$. 
\end{proof}

\begin{remark}
Under further assumptions on $X$, one can give more explicit  bounds on $\Delta$ in Theorem \ref{thm-BNmain}. For example, in addition to $\PP^2$ \cite{GouldLiu}, 
the Brill-Noether loci where $c_1({\bf v})$ is a minimal effective class have been studied for K3 surfaces by Leyenson \cite{Leyenson, Leyenson2} and Bayer \cite{BayerBN} and for abelian surfaces by Bayer and Li \cite{BayerLi} (see also \cite{CoskunNuer}).
\end{remark}

\bibliographystyle{plain}

\end{document}